\documentclass{ar-1col}
\usepackage[numbers]{natbib}
\usepackage{url}
\usepackage{amsfonts,amssymb}
\usepackage{bm}
\usepackage{amsmath}
\usepackage{amsthm}
\usepackage{mathrsfs}  
\usepackage{empheq}
\usepackage{cases}
\usepackage{dsfont}
\usepackage[inline]{enumitem}

\usepackage{algorithmicx}
\usepackage{algorithm}
\usepackage{algpseudocode}

\newcommand{\mR}{\mathbb{R}}

\newcommand{\spn}{\mathrm{span}} 
\newcommand{\im}{\mathrm{Im}} 
\newcommand{\diag}{\mathrm{Diag}}
\newcommand{\blkdiag}{\mathrm{blkdiag}}
\newcommand{\tr}{\mathrm{Tr}}

\newtheorem{thm}{Theorem}
\newtheorem{exm}{Example}

\newtheorem{defn}{Definition}

\jname{Annu. Rev. Control Robot. Auton. Syst.}
\jyear{2021}
\doi{10.1146/annurev-control-061820-083817}

\begin{document}
	
	\markboth{Cheng $\bullet$ Scherpen}{Model Reduction Methods for Network Systems}
	
	\title{Model Reduction Methods for Complex Network Systems}

	\author{X. Cheng$^1$, and J.M.A. Scherpen$^2$
		\affil{$^1$Control Systems Group, Department of Electrical Engineering, Eindhoven University of Technology, Eindhoven, the Netherlands; email: x.cheng@tue.nl}
		\affil{$^2$Jan C. Willems Center for Systems and Control, ENTEG, University of Groningen, the Netherlands; email: j.m.a.scherpen@rug.nl}}
	
	\begin{abstract} 
		Network systems consist of subsystems and their interconnections, and provide a powerful framework for analysis, modeling and control of complex systems. However, subsystems may have high-dimensional dynamics, and the amount and nature of interconnections may also be of high complexity. Therefore, it is relevant to study reduction methods for network systems. An overview on reduction methods for both the topological (interconnection) structure of the network and the dynamics of the nodes, while preserving structural properties of the network, and taking a control systems perspective, is provided. First topological complexity reduction methods based on graph clustering and aggregation are reviewed, producing a reduced-order network model. Second, reduction of the nodal dynamics is considered by using extensions of classical methods, while 
		preserving the stability and synchronization properties. Finally, a structure-preserving generalized balancing method for simplifying simultaneously the topological structure and the order of the nodal dynamics is treated. 
	\end{abstract}
	
	\maketitle
	

	\section{INTRODUCTION}
	
	The backbone of many modern technological systems is a network  system (system of systems), which bonds diverse multi-physics components together. Many large-scale systems can be modeled as  network systems which are composed of multiple subsystems interacting with each other via certain coupling protocols. Such systems are becoming ever more prevalent in various domains. Chemical reaction chains, cellular and metabolic networks, social networks, multi-robot coordination, and large-scale power grids
	are only a few examples \cite{mesbahi2010graph,Newman2003Review,ren2005survey,Rao2013graph,Ahsendorf2014GeneRegulation,proskurnikov2017tutorial,dorfler2012synchronization,Chow2013PowerReduction}. 
	
	However, with the increasing complexity of network scales and subsystem dynamics, the models describing the behavior of network systems can be of extremely high dimension. This will lead to serious scalability issues in simulation, optimization, transient analysis, and  control synthesis due to limited computational capability and storage capacity. These issues spur the development of methodologies on complexity reduction for large-scale network systems, aiming to acquire pertinent information of the systems in a computationally manageable fashion. 
	
	In the past few decades, a variety of theories and techniques for model reduction have been developed for generic dynamical systems, including Krylov-subspace methods (also known as moment matching), balanced truncation, and Hankel norm approximation, see    \cite{Antoulas05,Bai2002Krylov,Astolfi2010Moment,Moore81,glover1984all,besselink2013comparison} and the references therein. These conventional methodologies can generate reduced-order models that well approximate the input-output mapping of a high-dimensional system.  However, when addressing the model reduction problem of large-scale network system, we have to rethink about how to implement those methods in a structure-preserving manner. This is because analysis, control and monitoring of complex networks rely heavily on their interconnection structure \cite{summers2014optimal,gates2016control,Kim2018role,ishizaki2018graph}. Actually, preserving essential network configurations in the approximation of network systems presents the most challenging problem.  {Early work on controller reduction can be viewed as a predecessor of structure-preserving reduction of interconnected systems, which takes into account the coupling structure between plants and controllers \cite{Obinata2012MR,Mustafa1991Controller}. However, recent developments in large-scale networked systems have gone far beyond the simple closed-loop structure.}

	In this paper, we provide an overview of recent advances in dimension reduction of complex network systems. The complexity we consider consists of two aspects, namely, large-scale topology and high-dimensional subsystems (nodal dynamics), which lead to two types of model reduction problems in the context of network systems.
	The first one is focused on how to simplify a complicated network structure by reducing the number of nodes. Inspired by the classification and pattern recognition in data science and computer graphics \cite{SurveyClustering,Schaeffer2007SurveyClustering}, reduction methods based on \textit{clustering and aggregation} are mainstream for reducing the topological complexity. Most of the relevant work \cite{Schaft2014,Ishizaki2014,Monshizadeh2014,deng2014structure,Besselink2016Clustering,ChengTAC20172OROM,ChengTAC2018MAS,jongsma2018model,Martin2018scalefree} is treated in this paper. Besides, we also briefly review other topological methods, including the well-known \textit{singular perturbation approximation}. The second problem considers how to reduce the dimension of individual subsystems in a network. The relevant approximation 
	approaches for interconnected systems or coupled systems based on subsystem structuring have been of interest already for a long
	time \cite{Reis2008survey,Vandendorpe2008,sandberg2009interconnected}. Recent developments in \cite{Monshizadeh2013stability,ChengEJC2018Lure} further discuss diffusively coupled linear/nonlinear systems, where the reduction is performed on each subsystem in a way that certain properties of the entire network, such as synchronization and stability are retained. Furthermore, techniques in \cite{Ishizaki2016dissipative,ChengAuto2019} combine the complexity reduction of network structures and subsystem dynamics, also providing an attractive way to simplify the complexity of entire network systems.

	\section{FROM GRAPHS TO NETWORK SYSTEMS}
	In this section, we recapitulate some preliminaries on algebraic graph theory, and provide key concepts to model, analyze and design network systems. The graph-based modeling of network systems is then introduced. 
	We refer to e.g., \cite{Bullo2019Lectures,Wu2007Synchronization} for more details.  
	
	\subsection{Algebraic Graph Theory}
	The language of graphs is essential in modeling and control of networked systems, and it provides a natural tool for characterizing the interconnection structure of a network. 
	Any finite graph $\mathcal{G}$ can be featured by a finite and nonempty node set  $\mathcal{V}: = \{1, 2, ... , n\}$ and an edge set $\mathcal{E} \subseteq \mathcal{V} \times \mathcal{V}$. Depending on whether the edges have specific orientations, we have two basic categories of graphs:
	
	\subsubsection{Directed Graphs} Directed graph (in short, a \textit{digraph}) structures are found in various applications, 
	including biochemical reactions and social networks, see e.g.,\cite{Ahsendorf2014GeneRegulation,proskurnikov2017tutorial}, where the transmission of information or energy among network nodes is directional.
	This directionality can be encoded in the edges that are ordered pairs of elements of $\mathcal{V}$, and we say that there is an edge directed from node $i$ to node $j$ if $(i,j) \in \mathcal{E}$. A digraph $\mathcal{G}$ is called \textit{simple}, if it does not contain self-loops (i.e., $\mathcal{E}$ does not contain edges of the form $(i,i)$, $\forall~i$), and there exists exactly one edge directed from $i$ to $j$ if $(i,j) \in \mathcal{E}$. In a simple digraph $\mathcal{G}$, a node $i_n$ is \textit{reachable} from another node $i_0$, if there is a \textit{directed path} from $i_0$ to $i_n$. Here, this path is defined as a sequence of edges of the form $(i_{k-1}, i_k)$, $k = 1,...,n$, which joins a sequence of distinct nodes $i_0, i_1, ..., i_n$.
	
	Next, the connectivity notions for a digraph $\mathcal{G}$ are presented.
	\begin{enumerate*} [label=(\roman*)]
		\item $\mathcal{G}$ is strongly connected if any two nodes are reachable from each other;
		\item $\mathcal{G}$ is quasi strongly connected if all the nodes are reachable from a common node;
		\item $\mathcal{G}$ is weakly connected if its undirected version $\mathcal{G}_u: = (\mathcal{V}, \mathcal{E}_u)$ is strongly connected, where the set $\mathcal{E}_u \supseteq \mathcal{E}$ include both $(i,j)$ and $(j,i)$, if there is an edge $(i,j) \in \mathcal{E}$.
	\end{enumerate*}
	Note that any simple digraph $\mathcal{G}$ can be decomposed into a unique set of maximal strongly connected components (SCCs), which are the largest strongly connected subgraphs of $\mathcal{G}$. If an SCC has only outflows to other SCCs, it is called a root SCC  (RSCC). A weakly connected digraph may contain multiple RSCCs, while a quasi strongly connected digraph has only one RSCC.

	There are three matrices commonly used to characterize the topology of a digraph. The \textit{incidence matrix} $\mathcal{B} \in \mathbb{R}^{n \times |\mathcal{E}|}$ of $\mathcal{G}$ is defined such that
	$[\mathcal{B}]_{ij} = 1$, if edge $(i,j)\in \mathcal{E}$;  
	$[\mathcal{B}]_{ij} =  -1$, if edge $(j,i)\in \mathcal{E}$; 
	and 
	$[\mathcal{B}]_{ij} =  0 $ otherwise, 
	where each column indicates a directed edge.
	While this edge is assigned a positive value (weight), i.e., $\mathcal{G}$ is \textit{weighted}, we define a \textit{weighted adjacency matrix} $\mathcal{A}$, where $[\mathcal{A}]_{ij}$ is equal to the weight of the edge $(j,i)$ if $(j,i) \in \mathcal{E}$, and $[\mathcal{A}]_{ij}=0$ otherwise. Moreover, the weighted out-degree and in-degree matrices of $\mathcal{G}$ are the diagonal matrices defined by $D_{\mathrm{out}} : = \diag(\mathcal{A} \mathds{1})$ and $D_{\mathrm{in}} : = \diag(\mathds{1}^\top \mathcal{A} )$, respectively. A strongly connected digraph is called \textit{balanced}, if $D_{\mathrm{out}} = D_{\mathrm{in}}$. The \textit{Laplacian matrix} of a digraph $\mathcal{G}$ is defined as ${L}: = D_{\mathrm{out}} - \mathcal{A}$, and the elements of $L$ are given by
	\begin{equation} \label{defn:Laplacian}
	[{L}]_{ij} =  \begin{cases} 
	\sum_{j=1,j\ne i}^{n} \mathcal{A}_{ij}, & i = j\\
	-\mathcal{A}_{ij}, & \text{otherwise.}
	\end{cases} 
	\end{equation}
	Laplacian matrices are instrumental in  modeling various diffusion
	processes, e.g., \cite{fax2001graph,Mirzaev2013LaplacianDynamics,Rao2013graph}. A Laplacian matrix enjoys two fundamental properties: 
	\begin{enumerate*}[label=(\roman*)]
		\item $L \mathds{1} = 0$;
		\item $[L]_{ii} \geq 0$, $\forall~i \in \mathcal{V}$, and $[L]_{ij} \leq 0$, $\forall~i \ne j$.
	\end{enumerate*}
	Conversely, a real square matrix satisfying the two properties can also be interpreted as a Laplacian matrix that represents a weighted simple digraph. Note that Laplacian matrices are singular. If a weakly connected digraph has $m$ LSCCs, then its Laplacian matrix also has semisimple zero eigenvalues with multiplicity $m$, while all the other nonzero eigenvalues have positive real parts.

	\subsubsection{Undirected Graphs}
	Undirected graphs are commonly used to characterize interconnection structure of physical systems, e.g., power grids, RC circuits, and mass-damper systems,  \cite{dorfler2018electrical,Schaft2017modeling,Schaft2014,ChengTAC20172OROM}. An undirected graphs can be viewed as a special digraph, whose weighted adjacency matrix $\mathcal{A}$ (or Laplacian matrix $L$) is symmetric.  
	In this case, we can define a Laplacian matrix using an alternative formula:
	\begin{equation} \label{eq:LBWB}
	L = \mathcal{B} W \mathcal{B}^T,
	\end{equation}
	where $\mathcal{B}$ is the incidence matrix obtained by assigning an arbitrary orientation to each edge of $\mathcal{G}$, and $W := \diag(w_1, w_2, \cdots, w_{|\mathcal{E}|})$ with $w_k$ the weight associated to the edge $k$, for each $k = 1,2,..., |\mathcal{E}|$.
	If $\mathcal{G}$ is an undirected connected graph, the Laplacian matrix $L$ has the properties:
	\begin{enumerate*}[label=(\roman*)]
		\item ${L}^\top = {L}$ and $\ker L = \spn \{\mathds{1}\}   $  ;
		\item $[L]_{ij} \leq 0$ if $i \ne j$, and $[L]_{ii} > 0$.	
	\end{enumerate*}

	\subsection{Modeling of Network Systems}
	
	In the field of network science, the evolution of network topology over time is often specified as \textit{dynamics of networks}, \cite{Newman2006booknetwork}. Differently, for control systems, \textit{dynamics over networks} is of interest, where nodes represent individual dynamical systems that are coupled through edges \cite{Fagnani2017book}. We use the latter notion when referring to a \textit{network system}. Network system examples are chemical reaction networks, power grids, robotic networks, i.e., they  have a clear interconnection structure, physically or virtually. Additionally, network modeling is also applicable to spatially discretized systems that are originally described by PDE's, such as simple beam models or fluid dynamical systems.

	\subsubsection{Networks of Single-Integrators}

	The simplest network systems consider all the nodes being just single-integrators, namely, $\dot{x}_i(t) = v_i(t)$, where $x_i(t), v_i(t) \in \mathbb{R}$ are the state and input of node $i$. A digraph $\mathcal{G}$ then captures the interconnection topology of $n$ single-integrators, where the coupling rule is  
	\begin{equation}\label{eq:coupling}
	{v}_i(t) = -d_i x_i(t) +\sum_{j=1,j\ne i}^{n} [\mathcal{A}]_{ij}  x_j(t).
	\end{equation}
	In Equation~\ref{eq:coupling}, $d_i \in \mathbb{R}$ represents the state feedback gain, and $\mathcal{A}$ is the weighted adjacency matrix of $\mathcal{G}$, whose entry $[\mathcal{A}]_{ij}$ indicates the strength of the coupling between nodes $i$ and $j$. 
	Taking into account the external control signals $u(t)\in \mathbb{R}^p$ and measurements $y(t) \in\mathbb{R}^q$ of the network, we then derive a compact form for the network system as
	\begin{equation} \label{sys}
	\dot{x}(t) =  \Gamma x(t) + Fu(t),  \quad
	y(t) = Hx(t),
	\end{equation}
	where $\Gamma: = \mathcal{A} - D$ with $D  = \diag(d_1, ..., d_n)$, and $F\in \mathbb{R}^{n \times p}$, $H \in \mathbb{R}^{q \times n}$ are the input and output matrices, respectively. 
	
	Equation~\ref{sys} is regarded as a rather general representation for single-integrator networks, whose stability depends on the values in $D$ and $\mathcal{A}$.  If $D \geq 0$, $\Gamma$ becomes a \textit{Metzler matrix}, leading to the concept of monotone systems or positive systems \cite{rantzer2015positive,rantzer2018tutorial,ishizaki2015directed}.
	Particularly, if we choose $D > D_{\mathrm{out}}$ or $D > D_{\mathrm{in}}$, then the Metzler matrix $\Gamma$ is strictly row (column) diagonally dominant. Following the Gershgorin circle theorem \cite{Johnson1990Matrix}, $\Gamma$ is Hurwitz, leading to the asymptotic stability of the network system.
	If $D = D_{\mathrm{in}}$, or equivalently $\Gamma = - L^\top$ with $L$ the Laplacian matrix of $\mathcal{G}$, we have a network flow model \cite{Ahsendorf2014GeneRegulation,Mirzaev2013LaplacianDynamics}. Furthermore, if $D = D_{\mathrm{out}}$, i.e., $\Gamma = - L$, Equation~\ref{sys} becomes a \textit{consensus network}, or {continuous-time averaging systems} \cite{Bullo2019Lectures}. The coupling rule in Equation~\ref{eq:coupling} becomes 
	\begin{align} \label{eq:diffcoup}
	{v}_i(t) = - \sum_{j=1,j\ne i}^{n} [\mathcal{A}]_{ij} \left[x_i(t) - x_j(t)\right],
	\end{align} 
	which is known as the \textit{diffusive coupling} rule. For both $\Gamma = -L^\top$ and $\Gamma = -L$, the system in Equation~\ref{sys} is \textit{semistable} (or semi-convergent), i.e., $\lim_{t \rightarrow \infty}e^{\Gamma t}$ exists for any initial condition $x(0)$. Particularly, when $\mathcal{G}$ is strongly connected, then $\lim_{t \rightarrow \infty} e^{-L t} = \mathds{1} \omega^\top$, with $\omega$, satisfying $\mathds{1}^\top \omega = 1$, the left eigenvector of $L$ for eigenvalue $0$. 
	\begin{defn}
		A network system $\dot{x}(t) = \Gamma x (t)$ achieves  synchronization if 
		\begin{equation} \label{eq:consensus}
		\lim\limits_{t \rightarrow \infty}  \left[x_i(t) - x_j(t)\right]=0, \ \forall~i,j \in \mathcal{V},
		\end{equation}
		holds for all initial condition $x(0)$.
	\end{defn}
	
\begin{figure}[t]
		\centering
		\includegraphics[width=4in]{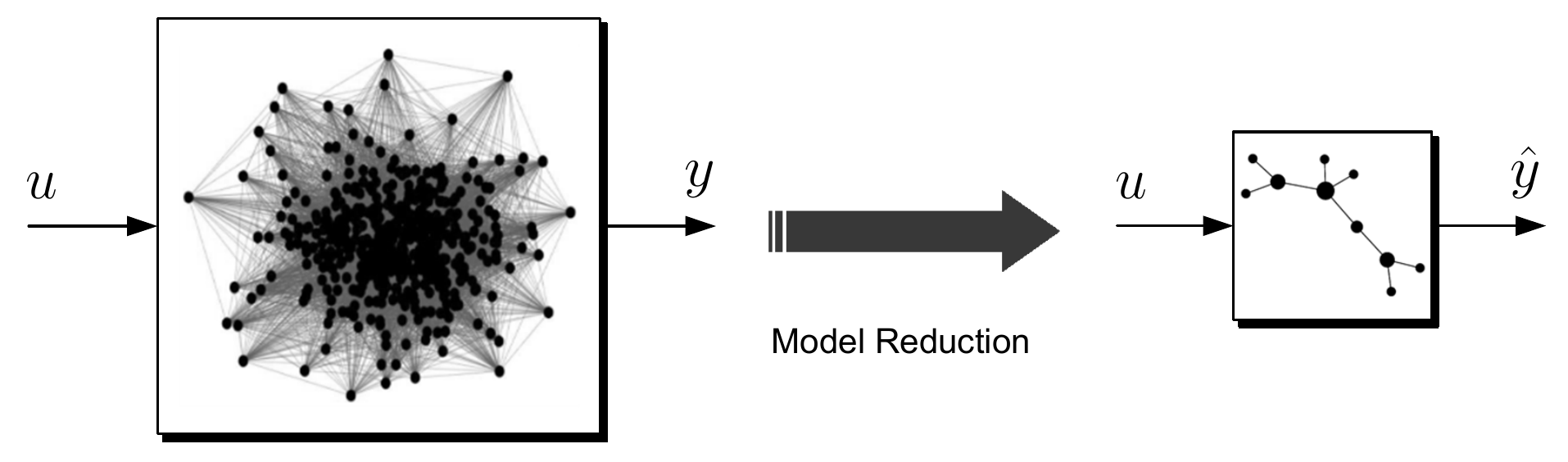}   
		\caption{Model reduction with preservation of network structure}
		\label{fig:networkreduction}
	\end{figure}

	
	The approximation of the network system~\ref{sys} aims for a reduced network consisting of a fewer number of nodes that captures essential properties of the original network. Specifically, a model reduction problem (\textbf{Figure~\ref{fig:networkreduction}}) is formulated to find a reduced-order model 
	\begin{equation} \label{sysr}
	\dot{\hat{x}}(t) = \hat{\Gamma} \hat{x}(t) + \hat{F} u (t), \quad \hat{y}(t) = \hat{H} \hat{x}(t), 
	\end{equation}
	where $\hat{x} \in \mathbb{R}^{r}$ ($r < n$), $\hat{y} \in \mathbb{R}^q$ such that \begin{enumerate*}[label=(\roman*)]
		\item $\hat{\Gamma} \in \mathbb{R}^{r \times r}$ is interpretable as a reduced graph, and
		\item the approximation error is minimized between the original and the reduced-order models.
	\end{enumerate*}
	The approximation error is usually evaluated by the $\mathcal{H}_\infty$ or $\mathcal{H}_2$ norms of 
	$\eta(s) - \hat{\eta}(s)$,
	\begin{equation} \label{eq:TFs}
	\eta(s): = H (sI_n - \Gamma)^{-1} F, \quad \hat{\eta}(s): = \hat{H} (sI_r - \hat{\Gamma})^{-1} \hat{F}.
	\end{equation}
	
	\subsubsection{Networked Linear Systems}
	
\begin{figure}[t]
		\centering
		\includegraphics[width=2.6in]{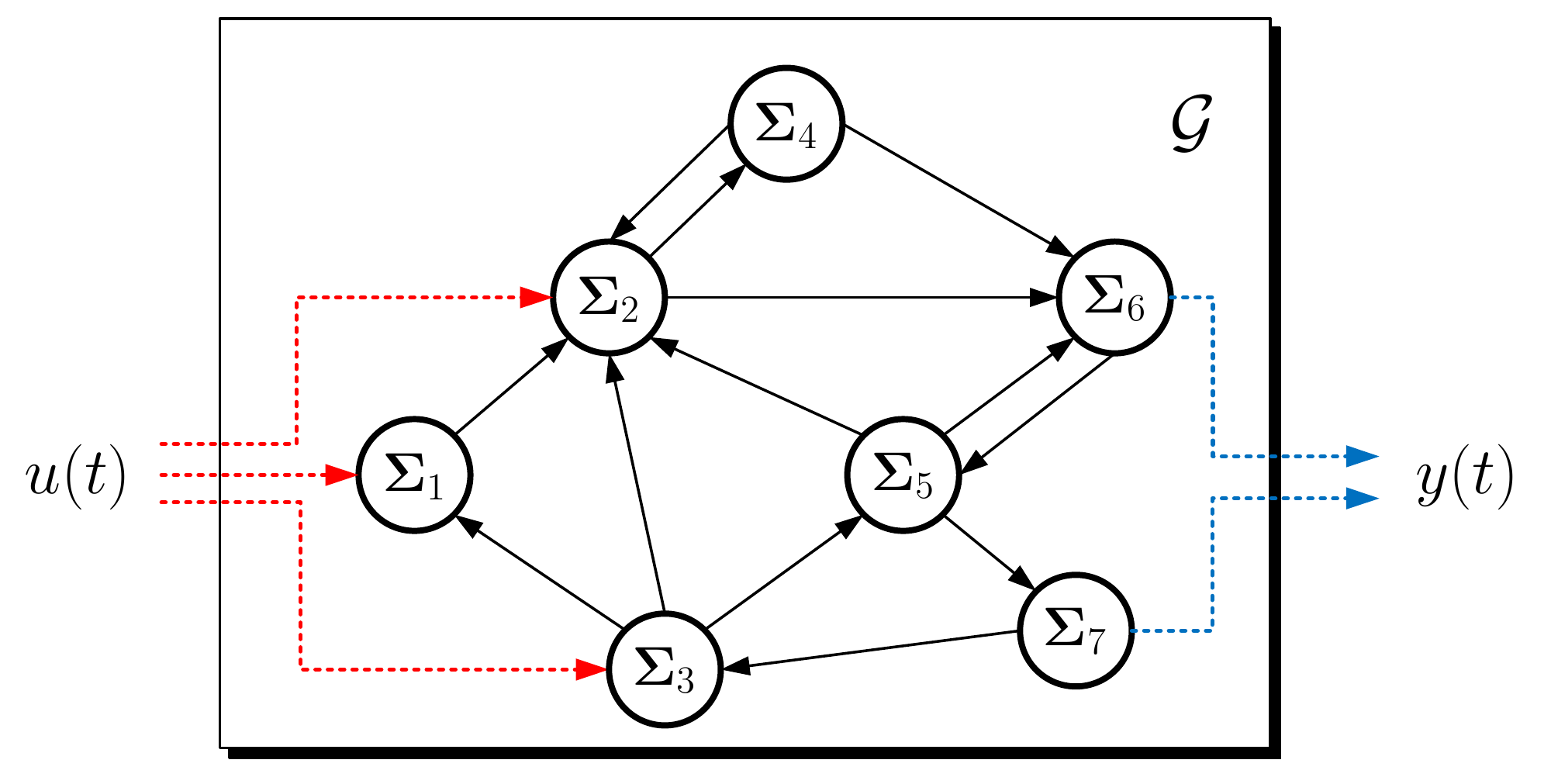}   
		\caption{An example of networked linear systems}
		\label{fig:interconnectedsys}
	\end{figure}

	The network model in Equation~\ref{sys} can be extended beyond single-integrators to consider each node as a high-order linear subsystem as 
	\begin{equation} \label{hsubsystem}
	\dot{x}_i(t) = A_i x_i(t) + B_i v_i(t), \quad
	y_i(t) = C_i x_i(t),
	\end{equation}
	where $x_i \in \mathbb{R}^{\ell_i}$, $v_i \in \mathbb{R}^{m_i}$, and $y_i \in \mathbb{R}^{\mu_i}$ are internal states, inputs and outputs, respectively. Suppose that $n$ subsystems are interconnected through the relations:
	$
	v_i(t) = \sum_{j=1}^{n} K_{ij} y_j(t) + F_i u(t), \quad y(t) = \sum_{i=1}^{n} H_i y_i(t),
	$ 
	with $K_{ij} \in \mathbb{R}^{m_i \times \mu_j}$ the coupling coefficient between nodes $i$ and $j$, where $K_{ij} = 0$ if and only if there are no signals passing from $j$ to $i$. The vectors $u(t)$ and $y(t)$ are denoted as external inputs and outputs. Combining this with Equation~\ref{hsubsystem}, we obtain a compact representation of the overall network system \cite{Reis2008survey,sandberg2009interconnected}:
	\begin{equation} \label{syshet}
	\dot{x}(t) = (A_{\mathrm{n}} + B_{\mathrm{n}} \Gamma C_{\mathrm{n}}) x(t) + B_{\mathrm{n}} F u(t),  \quad y(t) = H C_{\mathrm{n}} x(t),
	\end{equation}
	where $A_{\mathrm{n}}: = \blkdiag(A_1,...,A_n)$, $B_{\mathrm{n}}: = \blkdiag(B_1,...,B_n)$, $C_{\mathrm{n}}: = \blkdiag(C_1,...,C_n)$, and 
	\begin{equation*}
	\Gamma = \begin{bmatrix}
	K_{11} & \cdots & K_{1n} \\
	\vdots & \ddots & \vdots \\
	K_{n1} & \cdots & K_{nn} \\
	\end{bmatrix}, \quad
	F = \begin{bmatrix}
	F_1 \\ \vdots \\ F_n
	\end{bmatrix}, \quad
	H = \begin{bmatrix}
	H_1 & \cdots & H_n
	\end{bmatrix}.
	\end{equation*} 
	An example of a networked linear system containing six subsystems is shown in \textbf{Figure~\ref{fig:interconnectedsys}}. Networked linear systems of the form
	Equation~\ref{syshet} are also known as \textit{interconnected} or \textit{coupled} systems \cite{Reis2008survey,sandberg2009interconnected}. The subsystems in Equation~\ref{hsubsystem} could have different dynamics, in which sense the network is called \textit{heterogeneous}.  
	
	Homogeneous networks are defined when the dynamics of each node is identical:
	\begin{equation} \label{sysagent}
	\dot{x}_i(t) = A x_i(t) + B v_i(t), \quad
	y_i(t) = C x_i(t),
	\end{equation}
	where $x_i \in \mathbb{R}^{\ell}$, $v_i \in \mathbb{R}^{{m}}$ and $y_i \in \mathbb{R}^{{m}}$ are the internal state, input and output of node $i$, respectively. 
	Under a simple static output feedback interconnection similar to Equation~\ref{eq:coupling}, the dynamics of networked homogeneous linear systems is presented in compact form as
	\begin{equation} \label{sysh}
	\bm{\Sigma}:  
	\left \{
	\begin{array}{l}
	\dot{x}(t) = \left(I \otimes A + \Gamma  \otimes BC\right)x(t) + (F \otimes B)u(t),
	\\
	y(t)=(H \otimes C)x(t),
	\end{array}
	\right.
	\end{equation}
	with joint state vector $x(t) \in \mathbb{R}^{n\ell }$, external control inputs $u(t) \in \mathbb{R}^{pm}$ and external outputs $y \in \mathbb{R}^{qm}$. The matrix $\Gamma \in \mathbb{R}^{n \times n}$ indicates how the subsystems are interconnected. 
	
	A commonly studied network system in the form of Equation~\ref{sysh} is diffusively-coupled linear systems, where $\Gamma = - L$ is the Laplacian matrix of the underlying graph. Thus, the coupling rule among the nodes becomes 
	$
	{v}_i(t) = - \sum_{j=1,j\ne i}^{n} [\mathcal{A}]_{ij} \left[y_i(t) - y_j(t)\right].
	$  
	In this setting, the synchronization problem of the system $\bm{\Sigma}$ has been intensively studied in the literature, see e.g., \cite{ren2005survey,scardovi2008synchronization,consensus2010li,Wu2007Synchronization,Bullo2019Lectures}. 
	\begin{thm}  
		\label{thm:ConsCond}
		Consider a network of diffusively-coupled homogeneous linear systems described by Equation~\ref{sysh} with the symmetric Laplacian $L$. Let $0 = \lambda_1 < \lambda_2 \leq \cdots \leq  \lambda_n$ denote the
		eigenvalues of $L$. Then, the network system \ref{sysh} achieves  synchronization if and only if
		$A-\lambda_{i} BC$ is Hurwitz for all $i \in \{2,3,...,n\}$.  
	\end{thm}
	
	A sufficient condition for synchronization is provided by assuming that the subsystem $A,B,C$ is passive, i.e., there exists a symmetric positive definite matrix $K > 0$ that verifies 
	\begin{equation} \label{eq:passive}
	A^\top K + K A \leq 0, \ \text{and} \ C^\top = B K.
	\end{equation}
	Passivity is a natural property of physical systems, including mechanical systems, electrical networks, and thermodynamical systems  \cite{JanWillems1976}. With passivity we obtain a synchronization condition that is independent from the spectrum of the graph Laplacian \cite{arcak2007passivity,scardovi2008synchronization,chopra2012output}.
	
	\begin{thm} \label{thm:conspassive}
		Let $\Gamma = - L$ represent any connected undirected graph or any strongly connected digraph. If the subsystem $(A, B, C)$ in Equation~\ref{sysagent} is passive and observable, then the network system \ref{sysh} achieves synchronization.
	\end{thm}
	
	The model complexity of networked linear systems comes from two aspects: the dimension of subsystems and the  topological scale of the network. The first reduction problem is thus to \textit{reduce each subsystem} (or a subset of them) by taking into account the coupling structure in order to approximate the entire network system. For the reduction of heterogeneous network system~\ref{syshet}, the objective is to construct a network model composed of reduced-order subsystems $(\hat{A}_i, \hat{B}_i, \hat{C}_i)$, yielding an approximation of the entire system with the same form as \ref{syshet}, 
	where $\hat A_{\mathrm{n}}: = \blkdiag(\hat A_1,...,\hat A_n)$, 
	$\hat B_{\mathrm{n}}: = \blkdiag(\hat B_1,...,\hat B_n)$, 
	$\hat C_{\mathrm{n}}: = \blkdiag(\hat C_1,...,\hat C_n)$. The matrices $\Gamma$, $H$, and $F$ remain the same as the original one. Homogeneous network systems can be reduced in a similar manner such that each original subsystem $(A,B,C)$ is replaced by a lower-order approximation $(\hat{A},\hat{B},\hat{C})$:
	\begin{equation}  \label{syshsubsysr}
	\left \{
	\begin{array}{rcl}
	\dot{\hat{x}}(t) & = & \left(I_n \otimes \hat{A} -  {\Gamma} \otimes \hat{B}\hat{C}\right)\hat{x}(t) + ( {F} \otimes \hat{B}) u(t), \\
	\hat{y}(t) & = &  ( {H} \otimes \hat{C}) \hat{x}(t).
	\end{array}
	\right.
	\end{equation}
	
	The second reduction problem is focused on a \textit{simplification of the graph structure}, as illustrated in \textbf{Figure~\ref{fig:networkreduction}}. A resulting reduced-order model for networked homogeneous linear systems is in the form 
	\begin{equation} \label{sysrh}  
	\left \{
	\begin{array}{rcl}
	\dot{z}(t) & = & \left(I_r \otimes A - \hat{\Gamma} \otimes BC\right)z(t) + (\hat{F} \otimes B) u(t), \\
	\hat{y}(t) & = &  (\hat{H} \otimes C) z(t),
	\end{array}
	\right.
	\end{equation}
	where $\hat{\Gamma} \in \mathbb{R}^{r \times r}$ represents a reduced graph consisting of $r$ nodes. Let $G(s)$ and $\hat{G}(s)$ be the transfer matrices of the models \ref{sysh} and \ref{sysrh}, respectively. The objective now is to minimize the reduction error $G(s) - \hat{G}(s)$ with respect to certain norms. 

	\section{REDUCTION OF TOPOLOGICAL STRUCTURES} 
	\label{sec:topology}
	
	A powerful paradigm for simplifying a large-scale network is graph clustering. Graph clustering is a process of dividing a set of nodes into nonempty and disjoint subsets, where nodes in each subset are
	considered related by some similarity measure. Depending on the field, different names are used, including \textit{community detection} in social networks, and \textit{classification} in data science. Furthermore, it is closely related to \textit{unsupervised learning} in pattern recognition systems \cite{SurveyClustering,Schaeffer2007SurveyClustering}. 
	Generally, well-established  clustering algorithms (such as hierarchical clustering, spectral clustering or K-means clustering) were developed for static graphs or measured date. In this section, we present a series of clustering-based model reduction techniques for dynamic networks.

	\subsection{Clustering-Based Projection}
	\label{sec:clustering}
	
	Consider an LTI system with triplet $(A, B, C)$. The Petrov-Galerkin framework, \cite{Antoulas05} projects the state-space onto a lower dimensional subspace, resulting in a reduced-order model $(V^\dagger A V, V^\dagger B, C V)$, where $V$ is full column rank representing the basis of the subspace, and $V^\dagger$ is a left inverse of $V$, i.e. $V^\dagger V = I$. Clearly, the choice of $V$ is essential for obtaining the reduced-order model. For structure-preserving model reduction of network systems, $V$ can be constructed by considering an aggregation of node states.
	
	\begin{defn} \label{defn:cluster}
		Consider a graph $\mathcal{G}$ with node set $|\mathcal{V}| = n$. Graph clustering of $\mathcal{G}$ is a process that divides $\mathcal{V}$ into $r$ nonempty and disjoint subsets, denoted by  $   \mathcal{C}_1,\mathcal{C}_2,...,\mathcal{C}_r$, where $\mathcal{C}_i$ is called a \textit{cluster} (or a \textit{cell} of $\mathcal{G}$). 
		The {characteristic matrix} of the clustering $\{  \mathcal{C}_1,\mathcal{C}_2,...,\mathcal{C}_r\}$ is a binary matrix $\Pi \in \mathbb{R}^{n \times r}$ with \vspace{-0.3cm}
		\begin{equation} \label{partition}
		\mbox{} \hspace{1.5cm}[\Pi]_{ij} := \left\{ \begin{array}{ll}
		1, & \text{if node $i \in \mathcal{C}_j$,}\\ 0, & \text{otherwise.}
		\end{array}
		\right.
		\end{equation}  
	\end{defn}
	Note that each row of $\Pi$ has exactly one nonzero element, indicating that each node is assigned to a unique cluster. The number of nonzero elements in each column is the cardinality of the corresponding cluster. Specifically,
	$
	\Pi \mathds{1}_r = \mathds{1}_n~\text{and}~\mathds{1}^\top_n \Pi=\left[|\mathcal{C}_1|, |\mathcal{C}_2|,...,|\mathcal{C}_r| \right].
	$ 
	For any given undirected graph Laplacian $L$, the matrix $\Pi^\top L \Pi$ is a Laplacian matrix representing an undirected graph of smaller size. This important property allows for structure-preserving model reduction of network systems using $\Pi$ for the Petrov-Galerkin projection. To construct a reduced-order network system with $r$ nodes as in Equation~\ref{sys} or Equation~\ref{sysh}, we first have to find a clustering that partitions the nodes of a network into $r$ clusters. 
	
	Consider the network system~\ref{sys}, which is assumed to be semistable, see the sidebar on linear semistable systems. Then, the projection matrix is defined as $V V^\dagger \in \mathbb{R}^{n \times n}$ 
	\begin{equation} \label{eq:projection}
		V = N \Pi \in \mathbb{R}^{n \times r}, \quad V^\dagger: = (\Pi^\top M N \Pi)^{-1} \Pi^\top M \in \mathbb{R}^{r \times n},
	\end{equation}
	where $M$ and $N$ are nonsingular diagonal weighting matrices \cite{ChengAuto2020Gramian}. A reduced-order model is thereby obtained in the form of Equation~\ref{sysr}, 
	where $\hat{\Gamma} = V^\dagger \Gamma V$, $\hat{F} = V^\dagger F$, and $\hat{H} = H V$. 
	\begin{thm} \cite{ChengAuto2020Gramian}
		Let $\eta(s)$ and $\hat{\eta}(s)$ be the transfer matrices of the original and reduced-order models in Equations~\ref{sys} and \ref{sysr}. Then, a bounded $\mathcal{H}_2$ reduction error is guaranteed, i.e., $\eta(s) - \hat{\eta}(s) \in \mathcal{H}_2$, if there exist diagonal and positive definite matrices $M$ and $N$ such that 
		\begin{equation} \label{eq:H2condition}
		(\mathbf{e}_{i}-\mathbf{e}_{j})^\top N^{-1} \mathcal{J} = 0, \ \text{and} \    \mathcal{J} M^{-1} (\mathbf{e}_{i}-\mathbf{e}_{j}) = 0, 
		\end{equation}
		holds for each pair $i,j\in\mathcal{C}_k$, with any $k \in \{1,2,\cdots,r\}$.  
	\end{thm}
	
	\begin{textbox}[t]\section{Linear Semistable Systems and Pseudo Gramians}
		\label{sec:semistable}
		Semistability is a more general concept than asymptotic stability as it allows for multiple poles that are zero. The systems' trajectories thus may converge to a nonzero Lyapunov stable equilibrium \cite{bhat1999lyapunov,hui2009semistability}. Specifically, a linear system  $\dot{x}(t) = Ax(t)$ is \textit{semistable} if $\lim\limits_{t \rightarrow \infty} e^{At}$ is non-zero and exists for all initial states $x(0)$, or equivalently, the zero eigenvalues of $A$ are \textit{semisimple}, 
		and all the other eigenvalues have negative real parts. 
		
		It is well-known that the standard controllability and observability Gramians, \cite{Antoulas05}, are not well-defined for a semistable system. Therefore, in \cite{ChengAuto2020Gramian}
		the definition of pseudo Gramians is presented. 
		Consider a linear semistable system $(A,B,C)$.
		The \textit{pseudo controllability and observability Gramians} are defined as 
		\begin{equation}
		\label{defn:PseudoGramians} 
		\mathcal{P} = \int_{0}^{\infty} 
		(e^{A t}-\mathcal{J}) BB^\top (e^{A^\top t}-\mathcal{J}^\top) \mathrm{d}t,
		\quad
		\mathcal{Q} = \int_{0}^{\infty} 
		(e^{A^\top t}-\mathcal{J}^\top) C^\top C (e^{A t}-\mathcal{J}) \mathrm{d}t,
		\end{equation}
		respectively, where $\mathcal{J} : = \lim\limits_{t \rightarrow \infty} e^{At}$ is a constant matrix. The pseudo Gramians $\mathcal{P}$ and $\mathcal{Q}$ in Equation~\ref{defn:PseudoGramians} are well-defined for semistable systems. 
		The pseudo Gramians can be computed as $\mathcal{P} = \tilde{\mathcal{P}} - \mathcal{J} \tilde{\mathcal{P}} \mathcal{J}^\top$ and $\mathcal{Q} = \tilde{\mathcal{Q}} - \mathcal{J}^\top \tilde{\mathcal{Q}} {\mathcal{J}}$, where $\tilde{\mathcal{P}}$ and $\tilde{\mathcal{Q}}$ are arbitrary symmetric solution of the Lyapunov equations
		\begin{align}
		A \tilde{\mathcal{P}} + \tilde{\mathcal{P}} A^\top + (I-\mathcal{J})BB^\top(I-\mathcal{J}^\top) = 0,  \quad
		A^\top \tilde{\mathcal{Q}} + \tilde{\mathcal{Q}} A  + (I-\mathcal{J}^\top)C^\top C(I-\mathcal{J}) = 0,
		\end{align}
		respectively. The pseudo Gramians are useful for computing the $\mathcal{H}_2$ norm of a semistable system. 
		The transfer matrix $G(s)  \in \mathcal{H}_2$ if and only if $
		C \mathcal{J} B = 0
		$. 
		Furthermore, 
		$
		\lVert G(s) \rVert_{\mathcal{H}_2}^2 = \tr (C \mathcal{P} C^\top) = \tr (B^\top \mathcal{Q} B)
		$. 
	\end{textbox}
	
	If $\Gamma$ is Hurwitz, or $\Gamma = - L$ with $L$ the Laplacian of a connected undirected graph, we simply choose $M = N = I_n$, which always guarantees $\eta(s) - \hat{\eta}(s) \in \mathcal{H}_2$. A Hurwitz $\Gamma$ implies $\mathcal{J} = 0$, while a connected undirected graph yields $\mathcal{J} = \frac{1}{n} \mathds{1}_n \mathds{1}_n^\top$. 
	In \cite{ishizaki2015directed,ChengCDC2017Digraph,ChengTAC2020Weight}  dynamic networks having a strongly connected topology are treated, i.e.,  $\Gamma$ is irreducible, and has only one zero eigenvalue with corresponding left and right eigenvectors $\mu_l$ and $\mu_r$ 
	the so-called \textit{Frobenius eigenvectors} with all real and positive entries \cite{farina2011positive}.
	In this case $M = \diag(\mu_l)$ and $N = \diag(\mu_r)$ satisfy Equation~\ref{eq:H2condition} for any clustering. Furthermore, we have $\eta(s) - \hat{\eta}(s) \in \mathcal{H}_2$ and $N \Gamma^\top M  + M  \Gamma N \leq 0$. Following \cite{ChengAuto2020Gramian}, we can obtain a \textit{posteriori} bound on the reduction error as
	$
	\lVert \eta(s) - \hat{\eta}(s) \rVert_{\mathcal{H}_2} \leq 
	\gamma_s \sqrt{\tr(I-V V^\dagger) \mathcal{P}(I-V V^\dagger)^\top},
	$
	where $\mathcal{P}$ is the pseudo controllability Gramian of  system~\ref{sys}, and $\gamma_s  \in  \mathbb{R}_{>0}$ satisfies
	\begin{equation} \label{eq:LMI}
	\begin{bmatrix}
	N \Gamma^\top M  + M \Gamma N & M \Gamma  & (I - \mathcal{J}^\top) H^\top \\
	\Gamma^\top M & -\gamma_s I & H^\top\\
	H (I - \mathcal{J}) & H & -\gamma_s I
	\end{bmatrix} \leq 0. 
	\end{equation}
	There is a balanced graph representation of the digraph system \ref{sys} as follows:
	\begin{equation} \label{sysGbal}
	M N \dot{\xi}(t) = L_b \xi(t) + M F u(t), \quad y(t) = H N \xi(t),
	\end{equation}
	where $L_b: = M \Gamma N$ is the Laplacian matrix of the balanced digraph, and the resulting reduced-order model in Equation~\ref{sysr} becomes
	\begin{equation*}
	\Pi^\top M N \Pi \dot{\hat{\xi}}(t) = \Pi^\top L_b \Pi \hat{\xi}(t) + \Pi^\top M F u(t), \quad \hat{y}(t) = H N \Pi \hat{\xi}(t),
	\end{equation*}
	with $\Pi^\top L_b \Pi$ representing a reduced balanced digraph. 
	In \cite{ChengTAC2019Digraph}, a generalized balanced digraph is defined as a  weakly connected digraph in which each RSCC is balanced while removing all the non-RSCC nodes resulting in a generalized balanced graph representation similar to Equation~\ref{sysGbal}. 
	For networks with a weakly connected topology the error system generally $\eta(s) - \hat{\eta}(s) \notin \mathcal{H}_2$. Then, \textit{clusterability} is defined between two nodes $i,j$ if they satisfy Equation~\ref{eq:H2condition}. Clusterability of all nodes in each cluster then guarantees the stability of the error $\eta(s) - \hat{\eta}(s)$ \cite{ChengTAC2019Digraph}. 	
	
	
	\begin{exm} \label{exm1}
		\begin{figure}[t]
			\centering
			\includegraphics[width=4in]{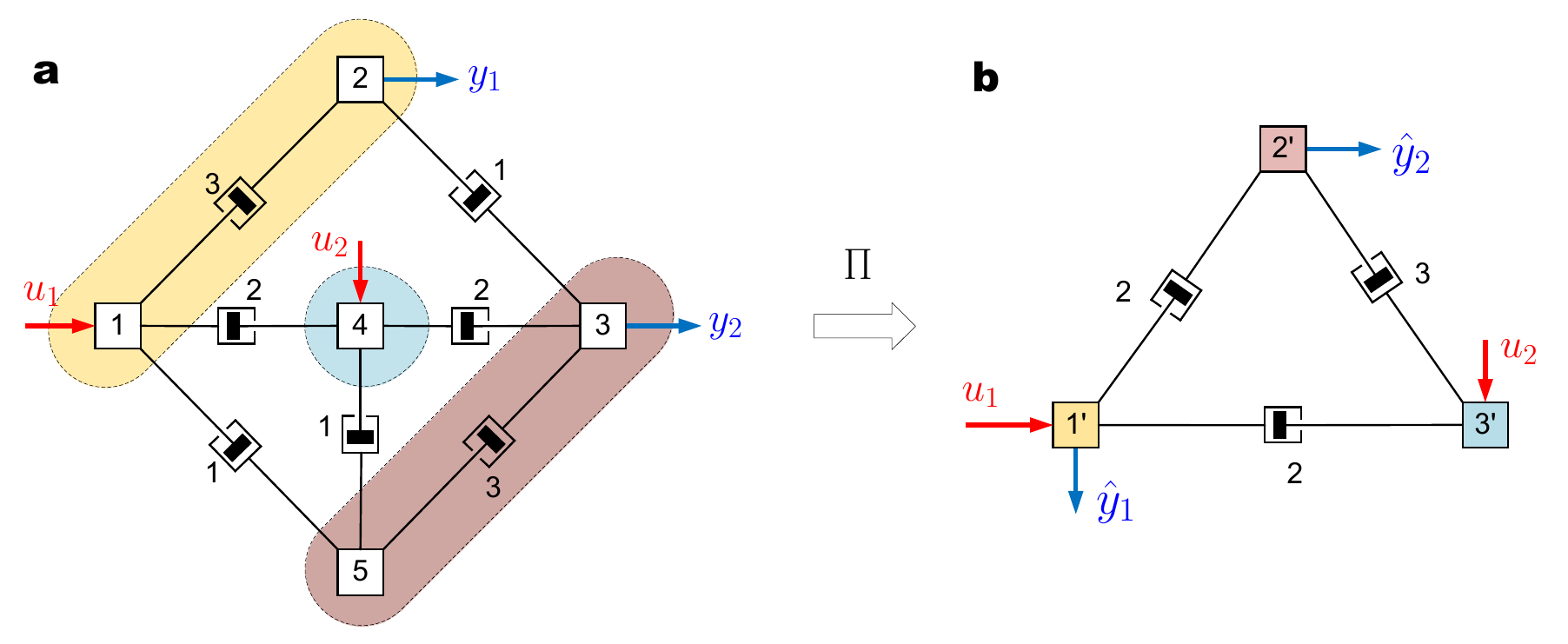}     
			\caption{An illustrative example of clustering-based model reduction of a mass-damper network system. (\textit{a}) The original network are divided into three clusters. (\textit{b}) The network representation of the reduced-order model.}
			\label{fig:massdamper}
		\end{figure}
		
		Consider the mass-damper system in  \textbf{Figure~\ref{fig:massdamper}a}, where the masses are interconnected via linear dampers. ${u}_1(t)$, $u_2(t)$ represent external forces, and $y_1(t)$ and $y_2(t)$ are measured velocities. Suppose that all the masses are identical, the network system in the form of Equation~\ref{sys} is obtained as   
		\begin{align*}
		\begin{bmatrix}
		\dot{x}_1 \\ \dot{x}_2 \\\dot{x}_3 \\\dot{x}_4 \\ \dot{x}_5
		\end{bmatrix} 
		= 
		\underbrace{-\begin{bmatrix}
			6  & -3&  0 & -2 & -1\\
			-3 &  4&  -1&  0 & 0\\
			0  &  -1&  6&  -2 & -3\\
			-2 &  0&  -2&  5 & -1\\
			-1 &  0&  -3&  -1 & 5\\
			\end{bmatrix}}_{\Gamma} 
		\begin{bmatrix}
		{x}_1 \\ {x}_2 \\{x}_3  \\{x}_4 \\ {x}_5 
		\end{bmatrix}
		+
		\underbrace{\begin{bmatrix}
			1 & 0 \\   
			0 & 0 \\
			0 & 0 \\
			0 & 1 \\
			0 & 0 \\
			\end{bmatrix}}_{F} 
		\begin{bmatrix}
		u_1  \\ u_2 
		\end{bmatrix},
		\quad
		\begin{bmatrix}
		y_1  \\ y_2 
		\end{bmatrix}  
		=
		\underbrace{\begin{bmatrix}
			0 & 1 & 0 & 0 & 0 \\
			0 & 0 & 1 & 0 & 0 
			\end{bmatrix}}_{H} 
		\begin{bmatrix}
		{x}_1  \\ {x}_2  \\{x}_3  \\{x}_4  \\ {x}_5 
		\end{bmatrix},
		\end{align*}
		where $-\Gamma$ is an undirected graph Laplacian, and the off-diagonal entry $[\Gamma]_{ij}$ represents the damping coefficient of the edge $(i,j)$. Consider $\{\mathcal{C}_1, \mathcal{C}_2, \mathcal{C}_3\} = \left\{ \{1,2\}, \{3,5\},\{4\}\right\}$ to be the clustering of the graph, which leads to the following characteristic matrix
		\begin{equation*}
		\Pi=\left[                 
		\begin{matrix} 
		1 & 1 &  0 & 0 & 0\\
		0 & 0 &  1 & 0 & 1\\
		0 & 0 &  0 & 1 & 0\\
		\end{matrix}
		\right]^\top.   
		\end{equation*}
{Therefore, a reduced-order network model is obtained as
\begin{align*} 
\underbrace{\begin{bmatrix}
	2 & 0 & 0\\
	0 & 2 & 0\\
	0  & 0 & 1\\
	\end{bmatrix}}_{\Pi^\top \Pi}
\begin{bmatrix}
\dot{z}_1 \\ \dot{z}_2 \\\dot{z}_3
\end{bmatrix}
= 
-\underbrace{\begin{bmatrix}
	4 &  -2  &   -2\\
	-2  &   5   & -3\\
	-2 &   -3  & 5
	\end{bmatrix}}_{\Pi^\top \Gamma \Pi} 
\begin{bmatrix}
{z}_1 \\ {z}_2 \\{z}_3 
\end{bmatrix}
+
\underbrace{\begin{bmatrix}
	1 & 0  \\
	0 & 1 \\
	0 & 1
	\end{bmatrix}}_{ \Pi^\top F} 
\begin{bmatrix}
u_1  \\ u_2 
\end{bmatrix},
\
\begin{bmatrix}
\hat{y}_1  \\ \hat{y}_2 
\end{bmatrix}  
=
\underbrace{\begin{bmatrix}
	1 & 0 & 0  \\
	0 & 0 & 1 
	\end{bmatrix}}_{H \Pi} 
\begin{bmatrix}
{z}_1  \\ {z}_2  \\{z}_3 
\end{bmatrix},
\end{align*}
where $-\Pi^\top \Gamma \Pi$ is again an undirected graph Laplacian.  	
To bring it in the form of Equation~\ref{sysr}, we can define $\hat{\Gamma}: = (\Pi^\top \Pi)^{-1} \Pi^\top \Gamma \Pi$, $\hat{F}: = (\Pi^\top \Pi)^{-1} \Pi^\top F$, and $\hat{H}: = H \Pi$. 
However, from $-\Pi^\top \Gamma \Pi$ it follows that this model allows for a physical interpretation, as shown in \textbf{Figure~\ref{fig:massdamper}b}: the nodes in each cluster are aggregated into a single node in the reduced network, while all edges connecting nodes from two distinct clusters are merged to a single edge linking the corresponding nodes in the reduced network.}
	\end{exm}
	
	Analogously, and beyond the single integrator case, a reduced-order model of networked homogeneous linear systems in Equation~\ref{sysh} can be formed using the Petrov-Galerkin projection framework of Equation~\ref{eq:projection}, which gives a reduced-order model in the form of Equation~\ref{sysrh},
	where $\hat{\Gamma} := \Pi^\top {L} \Pi $, $\hat{F}=\Pi^\top F$ and $\hat{H} = H\Pi$ are the same as in Equation~\ref{sysr}. The new state vector $z^\top(t): =  \left[z_1^\top(t)\ z_2^\top(t)\ ...\ z_r^\top(t) \right] \in \mathbb{R}^{r \ell}$, $z_i(t) \in \mathbb{R}^{\ell}$, $i=1,\ldots r$ represents an estimate of the state vector of the dynamics of all the nodes in the $i$-th cluster. Note that the extension of clustering-based approaches towards networks of  heterogeneous subsystems in Equation~\ref{syshet} remains an open problem. A major challenge lies in the representation of a cluster of nonidentical subsystems.

	Denote $G(s)$ and $\hat{G}(s)$
	as the transfer matrices of the models \ref{sysh} and \ref{sysrh}, respectively. The analysis of the reduction error $G(s) - \hat{G}(s)$ is more complicated than in the single integrator case, and for general subsystems, the reduction error $G(s) - \hat{G}(s)$ may not be stable.
	However, there is a theoretical guarantee if the subsystem $(A, B, C)$ in Equation~\ref{sysagent} are \textit{observable} and \textit{passive}. With Theorem~\ref{thm:conspassive} it can be verified that if the original network is undirected, or strongly connected, the reduced-order network system in Equation~\ref{sysrh} achieves synchronization, and $G(s) - \hat{G}(s) \in \mathcal{H}_2$, for any clustering $\Pi$  \cite{Besselink2016Clustering,ChengTAC2018MAS}.

	In the framework of clustering-based projection, the approximation error $\| G(s) - \hat{G}(s) \|_{\mathcal{H}_2}$  only depends on the choice of graph clustering. Thus, the most crucial problem in this framework is how to determine clusters of nodes to minimize the approximation error. In the following, we review several specific cluster selection approaches.

	\subsubsection{Almost Equitable Partitions} \textit{Almost equitable partitions} provide a graph clustering where nodes in the same cluster are connected to  other clusters in a ``similar'' fashion.  
	\begin{defn}
		Consider a weighted undirected graph $\mathcal{G}$ with adjacency matrix $\mathcal{A}$. A clustering $\{\mathcal{C}_1, \mathcal{C}_2, ..., \mathcal{C}_r\}$ is called an almost equitable partition if for any indexes $\mu, \nu \in \{1,2,...,r\}$ with $\mu \ne \nu$, it holds that
		$
		\sum_{k \in \mathcal{C}_\nu}^{}w_{ik} = \sum_{k \in \mathcal{C}_\nu}^{}w_{jk}
		$, $\forall~i,j \in \mathcal{C}_\mu$, where $w_{ij}: = [\mathcal{A}]_{ij}$.
	\end{defn}
	\begin{figure}[t]
		\includegraphics[width=2in]{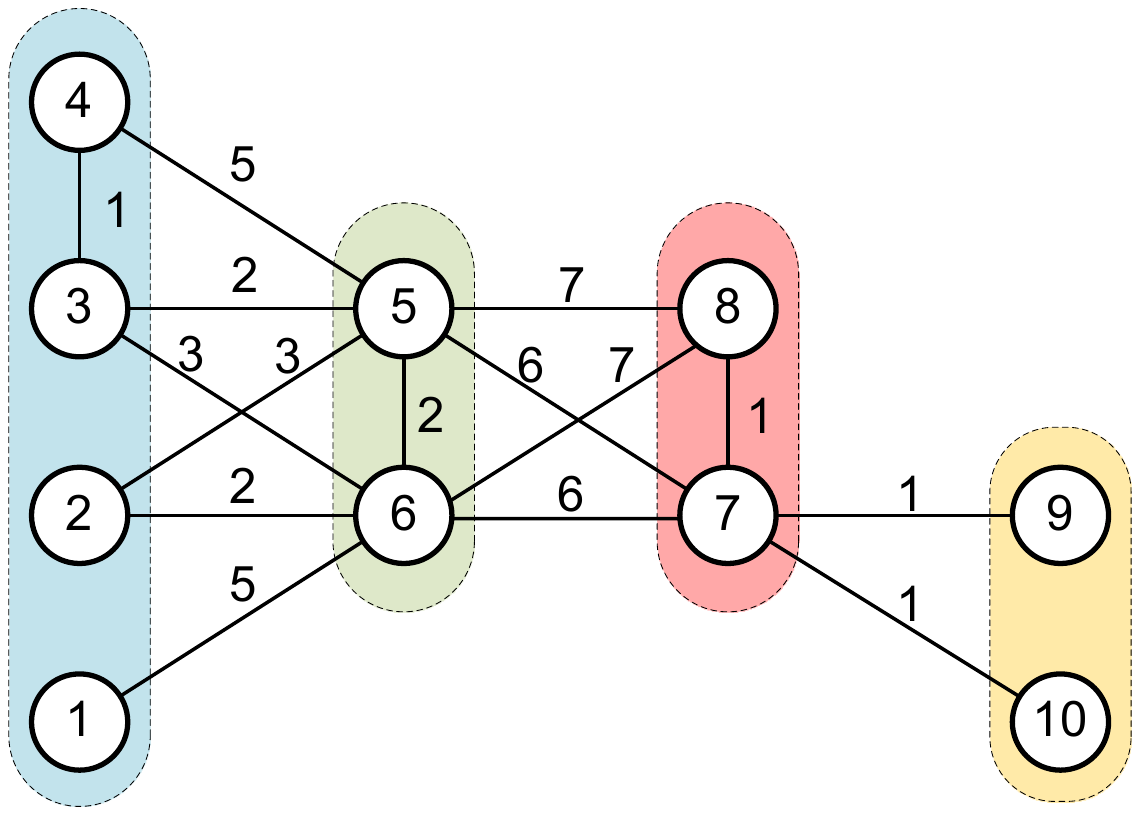}
		\caption{illustration of the almost equitable partition of a 10-node graph.}
		\label{fig:AEP}
	\end{figure}
	In \textbf{Figure~\ref{fig:AEP}}, an example of the almost equitable partition of an undirected graph is shown, \cite{Monshizadeh2014}. The nodes in a cluster have the same total edge weight to other clusters.  
	An almost equitable partition of an undirected graph has the key property that $\im(\Pi)$ is $L$-invariant,  i.e., 
	$
	L~\im(\Pi) \subseteq \im(\Pi),
	$ 
	where $L = L^\top$ is the Laplacian of an undirected graph, see e.g., \cite{zhang2013AEP,Monshizadeh2014}. Further, we have $L \Pi = \Pi \hat{L}$ with $\hat{L}: = (\Pi^\top \Pi)^{-1} \Pi^\top L \Pi$. A generalization of almost equitable partitions to digraphs is considered \cite{aguilar2017AEP}, where  nodes in the same cluster should have identical weighted out-degrees. Then, $\im(\Pi)$ is still $L$-invariant. 
	Now consider the following error system 
	\begin{equation} \label{eq:errorsystem}
	\Delta(s) = \eta(s) - \hat{\eta}(s) = \begin{bmatrix}
	H & -H \Pi
	\end{bmatrix}
	\begin{bmatrix}
	sI_n + L & 0 \\ 0 & sI_r + \Pi^\dagger L \Pi
	\end{bmatrix}^{-1}
	\begin{bmatrix}
	F \\ \Pi^\dagger F
	\end{bmatrix},
	\end{equation}
	where $\Pi^\dagger = (\Pi^\top \Pi)^{-1} \Pi^\top$ and
	$\eta(s)$, $\hat{\eta}(s)$ are the transfer matrices of Equation~\ref{eq:TFs}.  From the $L$-invariance it is verified that $\hat{\eta}(-s)^\top \Delta(s) = 0$, and thus 
	$\|\Delta(s)\|_{\mathcal{H}_2} = \|\eta(s)\|_{\mathcal{H}_2} - \|\hat{\eta}(s)\|_{\mathcal{H}_2} $.

	Furthermore, for a special output of the system~\ref{sys} 
	explicit expressions for the  reduction error $\|\Delta(s)\|_{\mathcal{H}_2}$ and $\|\Delta(s)\|_{\mathcal{H}_\infty}$ are provided in  \cite{Monshizadeh2014,jongsma2018model}.
	Further discussion on model reduction of networked symmetric linear system based on almost equitable partitions can be found in \cite{jongsma2018model}.
	Although an almost equitable partition  as a particular clustering offers us analytical expression for the reduction error, it does not necessary lead to a small error. In fact, the methods in \cite{Petar2015CDC,ChengECC2019weight} provide significantly lower errors via alternative choices of clustering for some examples.  Moreover, how to find all almost equitable partitions for a large-scale graph is generally a rather difficult and computationally expensive problem \cite{Monshizadeh2014}.
	
	\subsubsection{Tree Networks}
\begin{figure}[h]
		\includegraphics[width=2in]{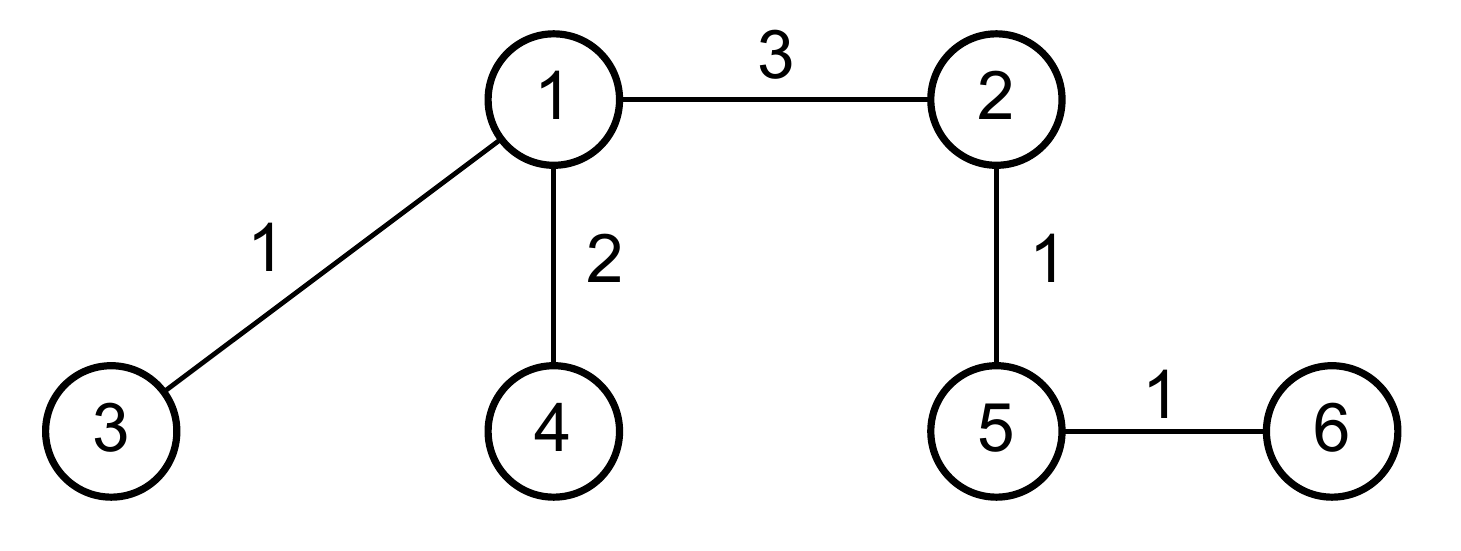}
		\caption{A tree with 6 nodes.}
		\label{fig:tree}
	\end{figure}

	Next, we focus on a particular class of undirected networks with tree topology. In graph theory, a \textit{tree} is a connected undirected graph in which there is only one path between any two nodes. An example of an undirected tree is shown in \textbf{Figure~\ref{fig:tree}}.
	
	Clearly, a tree $\mathcal{T}$ with $n$ nodes has exactly $n-1$ edges. Let $\mathcal{B} \in \mathbb{R}^{n \times (n-1)}$ be the incidence matrix of  $\mathcal{T}$. Relevant to the expression of the graph Laplacian $L$ in Equation~\ref{eq:LBWB}, we define an \textit{edge Laplacian} as $ 
	L_{\mathrm{e}} = \mathcal{B}^\top \mathcal{B} W \in \mathbb{R}^{(n-1) \times (n-1)},
	$ 
	where $W$ is the diagonal edge weight matrix. Observe that $L_{\mathrm{e}}$ is full rank and has all eigenvalues real and positive. The eigenvalues of $L_{\mathrm{e}}$ coincide with the nonzero eigenvalues of $L$, the Laplacian matrix of $\mathcal{T}$. 
	
	Consider the Laplacian dynamics in Equation~\ref{sys} where $\Gamma = - L$
	is an undirected graph Laplacian in Equation~\ref{eq:LBWB}. 
	Applying the transformation $x_{\mathrm{e}} = \mathcal{B}^\top x$ then leads to the so-called \textit{edge agreement protocol} \cite{mesbahi2010graph,zelazo2013performance}: 
	$
	\dot{x}_{\mathrm{e}}(t) = - L_{\mathrm{e}} x_{\mathrm{e}}(t) + \mathcal{B}^\top F u(t),
	$ 
	which is asymptotically stable and minimal. Network reduction approaches can be developed based on edge operations. For example, \cite{leiter2017graph} provides a greedy algorithm for edge-based contraction to simplify the graph topology. 
	A more general form of the edge agreement protocol is derived when subsystems are taken into account.
	In \cite{Besselink2016Clustering} the edge system of a network system with $\Gamma = -L$ representing an undirected tree graph is defined as
	\begin{equation} \label{sys:edge}
	\bm{\Sigma}_\mathrm{e}: \begin{cases}
	\dot{x}_\mathrm{e}(t) = (I_{n-1} \otimes A - L_{\mathrm{e}} \otimes BC) x_\mathrm{e}(t) + (\mathcal{B}^\top F \otimes B) u(t), \\
	y_\mathrm{e}(t) = (H \mathcal{B} W \otimes C) x_\mathrm{e}(t),
	\end{cases}
	\end{equation}
	where $x_\mathrm{e} = (\mathcal{B}^\top \otimes I) x \in \mathbb{R}^{(n-1)\ell}$. 
	Assuming that the subsystem $(A, B, C)$ in Equation~\ref{sysagent} is \textit{passive} and \textit{minimal}, we have the synchronization property of the network system $\bm{\Sigma}$ from Theorem~\ref{thm:conspassive}. 
	It then follows that the edge system $\bm{\Sigma}_\mathrm{e}$ is asymptotically stable. Furthermore, we can define a pair of generalized controllability and observability Gramians of the edge system  $\bm{\Sigma}_\mathrm{e}$ as follows, \cite{Besselink2016Clustering}, 
	$
	P_\mathrm{e} : = X \otimes K^{-1}, \quad Q_\mathrm{e} : = Y \otimes K,
	$ 
	where $K > 0$ satisfies \ref{eq:passive} for the passive subsystem, and $X >0$ and $Y >0$ are solutions of the following inequalities:
	\begin{equation*}
	-L_{\mathrm{e}} X - X L_{\mathrm{e}}^\top + \mathcal{B}^\top FF^\top \mathcal{B} \leq 0, \quad
	-L_{\mathrm{e}}^\top Y - Y L_{\mathrm{e}} + W\mathcal{B}^\top H^\top H \mathcal{B}W \leq 0. 
	\end{equation*}
	The matrices $X$ and $Y$ admit a diagonal structure: 
	$
	X = \diag (\xi_1, \xi_2, ..., \xi_{n-1}), \quad Y = \diag (\eta_1, \eta_2, ..., \eta_{n-1}),
	$ 
	where the ordering $\xi_i \eta_i \geq \xi_{i+1} \eta_{i+1}$ is imposed. The value of $\xi_i \eta_i$ can be roughly viewed as an indication for the importance of the $i$-th edge, since similar to balanced realization theory $\xi_i$ and $\eta_i$ are related to controllability and observability properties of the edges. Reduction by truncation methods then are equivalent to aggregating nodes connected by the truncated edges, and moreover \cite{Besselink2016Clustering} provides an \textit{a priori} upper bound for the approximation error in terms of the $\mathcal{H}_\infty$ norm 
	$
	\| G(s) - \hat{G}(s) \|_{\mathcal{H}_\infty} \leq 
	2 \left( \sum_{i = r }^{n-1} [L_{\mathrm{e}}^{-1}]_{ii} \sqrt{\xi_i \eta_i}\right),
	$ 
	where $r$ is the number of nodes in the reduced network.
	
	For an extension beyond tree graphs a major challenge lies in the characterization of the edge system, \cite{Besselink2016Clustering}.
	
	\subsubsection{Dissimilarity-Based Clustering} 
	
	For generic network systems, we may resort to a dissimilarity-based clustering approach presented in e.g., \cite{ChengECC2016,ChengTAC2018MAS,ChengTAC20172OROM,ChengTAC2019Digraph}. In line with data classification or pattern recognition in the other domains, dissimilarity-based clustering for dynamic networks starts with a proper metric that quantifies the difference between any pair of nodes (subsystems) in a network. For static graphs, dissimilarity (or distance) between two nodes or data points can be computed in a vector space using some sort of metric,   \cite{SurveyClustering,Schaeffer2007SurveyClustering}. 
	Considering dynamic	networks with external inputs, the dissimilarity metric is then featured in a function space \cite{ChengTAC2018MAS,ChengTAC20172OROM,ChengTAC2019Digraph}.  
	\begin{defn} \label{defn:dissim}
	Consider the network system of Equation~\ref{sysh}, the {dissimilarity} between nodes $i$ and $j$ is defined as 
	\begin{equation} \label{eq:dissim}
	\mathcal{D}_{ij} = \mathcal{D}_{ji} := \lVert \eta_i(s) - \eta_j(s) \rVert_{\mathcal{H}_2},
	\end{equation}
	where
	$
	{\eta}_i(s) : = (\mathbf{e}_i^\top \otimes C) (sI_{n\ell} - I_n \otimes A  + \Gamma \otimes BC)^{-1}(F \otimes B).
	$ 
	Consider a single-integrator network in Equation~\ref{sys}, node dissimilarity can be simply defined as 
	\begin{equation} \label{eq:dissim1}
	\mathcal{D}_{ij} = \lVert (\mathbf{e}_i -\mathbf{e}_j)^\top (sI_n - \Gamma)^{-1} F \rVert_{\mathcal{H}_2}.
	\end{equation}
\end{defn}

The transfer matrix $\eta_i(s)$  maps the external control signal $u \in \mathbb{R}^{p}$ to the measured output of the $i$-th node $y_i \in \mathbb{R}^{m}$. Thus $\eta_i(s)$ can be interpreted as the behavior of the $i$-th node with respect to the external inputs, while the dissimilarity measure in Equation~\ref{eq:dissim} indicates how different two nodes behave. 
The location of inputs and network topology determine the value of dissimilarity. 
Dissimilarity can also be defined in terms of other function norms, e.g., the $\mathcal{H}_\infty$ norm. However, to compute the dissimilarity between every pair of nodes in a large-scale network, the  $\mathcal{H}_2$ norm of a stable LTI system can be characterized by its Gramians \cite{Antoulas05}, whereas for other norms there is no such characterization, making them computationally less feasible.  

Note that the dissimilarity in Equation~\ref{eq:dissim} or Equation~\ref{eq:dissim1} is only well defined when $\eta_i(s) - \eta_j(s) \in \mathcal{H}_2$. This condition is guaranteed for network systems that are asymptotically stable or achieve synchronization \cite{ChengTAC2018MAS}. 
For instance, if $\Gamma$ in Equation~\ref{sys} is Hurwitz, then  Equation~\ref{eq:dissim1} immediately becomes
$
\mathcal{D}_{ij}  =  \sqrt{(\mathbf{e}_i -\mathbf{e}_j)^\top P (\mathbf{e}_i -\mathbf{e}_j)},
$
where $P$ is the controllability Gramian of the network system~\ref{sys}, and $P$ is computed as the unique solution of the Lyapunov equation $\Gamma P + P \Gamma^\top + F F^\top = 0$. If system~\ref{sys} is semistable, we resort to pseudo Gramians defined by Equation \ref{defn:PseudoGramians} for the computation of the $\mathcal{H}_2$ norm.

For example, a single-integrator network in Equation~\ref{sys} is considered with $\Gamma = - L$ the Laplacian matrix of an undirected graph. Following \cite{ChengCDC2016Gramian,ChengAuto2020Gramian}, the pseudo controllability Gramian of the system is computed as $\mathcal{P} = \mathcal{J} \tilde{\mathcal{P}} \mathcal{J}^\top$, where $\tilde{\mathcal{P}}$ is a solution of
\begin{equation} \label{eq:J}
- L \tilde{\mathcal{P}} - \tilde{\mathcal{P}} L + (I - \mathcal{J})  F F^{\top}  (I - \mathcal{J}) = 0, \quad \mathcal{J}: = \frac{1}{n} \mathds{1} \mathds{1}^\top.
\end{equation} 
The dissimilarity in Equation~\ref{eq:dissim1} is thereby obtained as  
\begin{equation} \label{eq:Dij1}
\mathcal{D}_{ij} = \sqrt{(\mathbf{e}_i - \mathbf{e}_j)^\top \mathcal{P} (\mathbf{e}_i - \mathbf{e}_j)}.
\end{equation} 
Further consider the network system~\ref{sysh} with $\Gamma = - L$ a symmetric Laplacian matrix. Assume that the network achieves synchronization. In this case, 
\cite{ChengTAC2018MAS} provides the expression for Equation~\ref{eq:dissim} as 
\begin{equation} \label{eq:Dij2}
\mathcal{D}_{ij} = \sqrt{\tr(\Psi_{ij} \bar{\mathcal{P}} \Psi_{ij}^\top)},
\end{equation}
where $\Psi_{ij}$ is defined with help of the output matrix $C$ and $\bar{\mathcal{P}} \in \mathbb{R}^{(n-1) \ell \times (n-1)\ell}$ is the unique solution of a Lyapunov equation with matrices built from the system matrices. 

The dissimilarity in Definition~\ref{defn:dissim} is a  pairwise measure that shows how close the behavior of two subsystems is, thus taking dynamics into account. This is significantly different from conventional node dissimilarity in data science or computer graphics \cite{SurveyClustering,Schaeffer2007SurveyClustering}. Nevertheless, we can still follow similar clustering procedures or algorithms for data sets or static graphs.  

Formally, given a network, the goal of clustering is to divide the nodes into clusters such that the elements assigned to a particular cluster are similar in a predefined metric. However, clustering with respect to a distance metric is generally an NP-hard combinatorial optimization problem, which is commonly solved by approximation algorithms. Here, we review two of such algorithms and their adaption to the clustering of network systems.

Agglomerative \textit{hierarchical clustering} is a method that produces multi-level clusters. The key of this method is to define the \textit{proximity} between two {clusters} $\mathcal{C}_i$ and $\mathcal{C}_j$. There are several  alternatives for such definition, 
such as considering the minimum, maximum or average dissimilarity between any node in $\mathcal{C}_i$ and any node in $\mathcal{C}_j$, respectively. The proximity of two clusters allows us to identify a pair of clusters with the smallest proximity and merge them into a single cluster. This operation is executed iteratively to generate clusters in a hierarchy structure, which is visualized as a \textit{dendrogram}, see an example shown in \textbf{Figure~\ref{fig:dendrogram}}. 
This is a tree-like diagram that records the sequences of cluster merges.
A graph clustering is obtained by cutting the dendrogram at the desired
level, then each connected component forms a cluster. The implementation of hierarchical clustering to the model reduction of network systems can be found in e.g., \cite{ChengACOM2018Power,ChengTAC20172OROM,ChengCDC2017Digraph,kawano2019data}.

\begin{exm}
	
	Consider the networked mass-damper system in Example~\ref{exm1}. The dissimilarity matrix can be computed as in Equation~\ref{eq:Dij1}, which yields
	\begin{equation*}  
	\mathcal{D}=\left[                
	\begin{array}{rrrrr}   
	0  &  0.2494  &  0.3154  &  0.3919  &  0.4142\\
	0.2494 &        0  &  0.2119  &  0.3688  &  0.3842\\
	0.3154 &   0.2119  &       0  &  0.2410  &  0.2394\\
	0.3919 &   0.3688  &  0.2410   &      0  &  0.0396\\
	0.4142 &   0.3842  &  0.2394  &  0.0396   &      0
	\end{array}
	\right].           
	\end{equation*} 
	We use the average-link to define the cluster proximity, and a dendrogram is generated as depicted in \textbf{Figure~\ref{fig:dendrogram}}, showing how clusters are merged hierarchically.  The dashed line cuts the dendrogram at a chosen level such that three clusters are formed:
	$ \{ 1\},\{ 2, 3\},\{ 4,5\}$. 
	
	\begin{figure}[h]\centering
		\includegraphics[scale=.4]{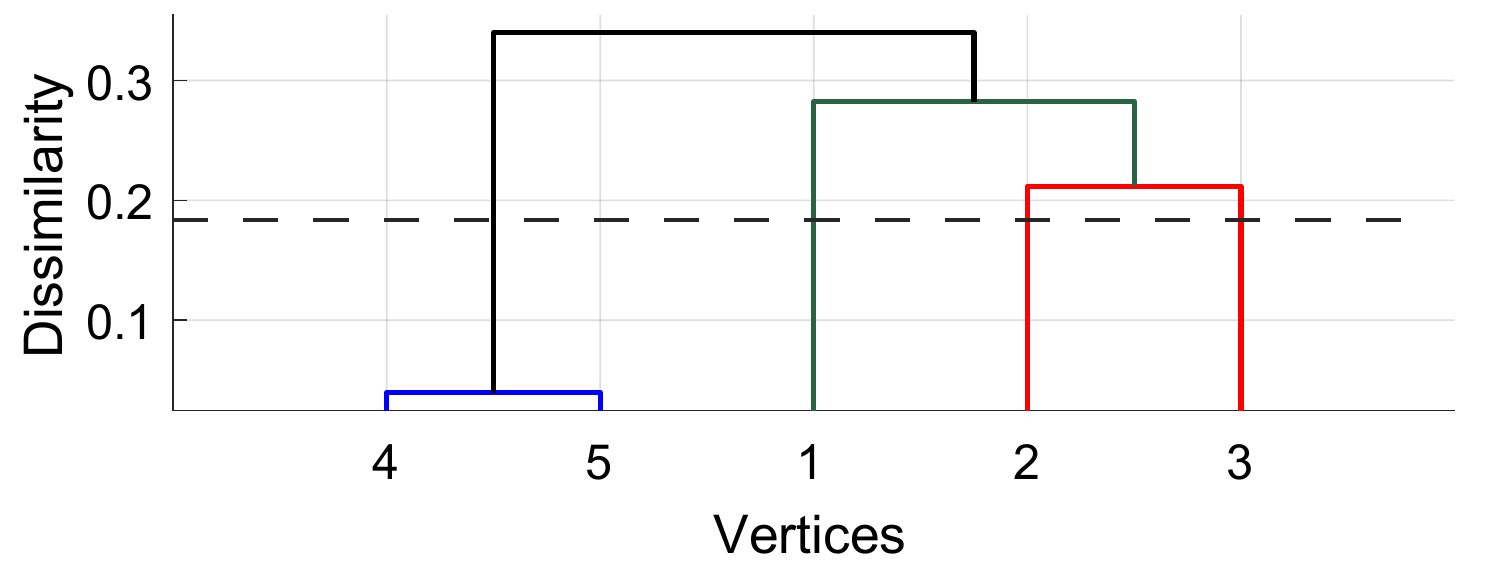}    
		\caption{Dendrogram illustrating the hierarchical clustering of the networked mass-damper system. The horizontal axis is
			labeled by node numberings, while the vertical axis represents
			the proximity of clusters. The level at which branches merge indicates the proximity between two clusters. }
		\label{fig:dendrogram}
	\end{figure}    
\end{exm}

\textit{K-means clustering} is a typical centroid-based partitioning method   \cite{SurveyClustering,zha2002spectral}, by which a cluster is constructed such that all the nodes within the
cluster is more similar to the \textit{centroid} of this cluster than to the centroid of any other clusters. For a network system of Equation~\ref{sysh}, the centroid of a cluster $\mathcal{C}_k$ can be defined as
\begin{equation} \label{eq:centroid}
\mu(\mathcal{C}_k) = \frac{1}{|\mathcal{C}_k|} \sum_{i \in \mathcal{C}_k}^{ } \eta_i(s),
\end{equation} 
with $\eta_i(s)$ defined in Equation~\ref{eq:dissim}. Given a network of $n$ nodes, K-means clustering aims to partition the nodes into $r$ subsets so as to minimize the following objective function:
$ 
\arg \min \sum_{k =1}^{r} \sum_{i \in \mathcal{C}_k}^{ } 
\| \eta_i(s) - \mu(\mathcal{C}_k) \|_{\mathcal{H}_2}^2,
$ 
in which $\eta_i(s) - \mu(\mathcal{C}_k) \in \mathcal{H}_2$ holds for synchronized networks. This problem can be solved using a simple iterative algorithm. First, take an initial $r$ clusters of a given network, and specify the centroid as in Equation~\ref{eq:centroid} for each cluster. Then, compute the dissimilarity between every node $i$ and the $r$ centroids, and assign node $i$ to the cluster whose centroid is the closest to $i$. Finally, we update the cluster centroids accordingly and repeat the steps until convergence.

Note that the formation of clusters in hierarchical clustering or K-means clustering solely relies on the dissimilarity measures and thus does not take into account the connectedness of nodes within a same cluster. It is worth noting that both methods can be modified to produce clusters of a graph, where each cluster forms a connected subgraph, \cite{ChengTAC2018MAS,ChengTAC2019Digraph,UmarCDC2019Clustering}. 

Generally, an upper bound on the reduction error $G(s) - \hat{G}(s)$, with $G(s)$ and $\hat{G}(s)$ the transfer matrices of systems~\ref{sysh} and \ref{sysrh}, is not easy to obtain. We thus impose extra assumptions on the network system~\ref{sysh}: $\Gamma = -L$ represents a connected undirected graph, and $A$ in Equation~\ref{sysagent} satisfies $A + A^\top < 0$. Then, a posteriori error bound is given as
$
\lVert G(s) - \hat{G}(s) \rVert_{\mathcal{H}_2} <  \gamma   \sum_{k=1}^{r} |\mathcal{C}_k| \cdot \max\limits_{i,j \in \mathcal{C}_k} \mathcal{D}_{ij},
$ 
where $\gamma  \in \mathbb{R}_{>0}$ only depends on the original system~\ref{sysh} and satisfies an LMI, 
\cite{Ishizaki2016dissipative,ChengTAC2018MAS}.

The most crucial part in dissimilarity-based clustering is to properly define the dissimilarity of nodes and clusters. For LTI network systems, dissimilarity can be defined using transfer matrices, which is applicable to different types of dynamical networks, see, e.g., \cite{ChengTAC20172OROM,ChengTAC2019Digraph,ChengTAC2019Digraph,ChengACOM2018Power} for more generalizations to second-order networks, directed networks and controlled power networks. However, how to extend the dissimilarity-based clustering to network systems containing nonlinearities still needs further exploration. One potential solution resorts to the \textit{dc gain} of monotone systems, that is introduced \cite{kawano2019data}, where the dc-gain can   
be regarded as indicator of the node importance. 

\subsection{Clustering Meets Optimization}
\label{sec:optimization}
In the previous section, we have reviewed how to select clusters and construct a reduced-order network model using the clustering-based projection. In this section, we formulate the model reduction problem from the perspective of  optimization, that is  to construct a lower-order network model which minimizes a certain  reduction error.

\subsubsection{Reducibility and an $\mathcal{H}_2$ Error Bound}
The pioneering work on clustering-based model reduction of dynamic networks in \cite{Ishizaki2014,ishizaki2015directed} introduces a notion of cluster reducibility, which is relevant to the classic notions \textit{exact aggregation} and \textit{approximate aggregation} from the control and model reduction literature \cite{aoki1968aggregation,feliachi1987interconnected}.

Consider the network system in Equation~\ref{sys} with $\Gamma$ a Hurwitz, Metzler and symmetric matrix and $F \in \mathbb{R}^n$. A cluster $\mathcal{C}_k$
is said to be \textit{reducible} if there exist a scalar rational function $g_\star(s)$ and a vector $p_k \in \mathbb{R}^{|\mathcal{C}_k|}$ such that
$
I({\mathcal{C}_k})   g(s) = p_k   g_\star(s),
$
where $g(s): = (sI - \Gamma)^{-1} F$, and $I({\mathcal{C}_k})$ denotes the matrix composed of the column vectors of $I_n$ compatible with the set $\mathcal{C}_k$. Reducibility reflects the uncontrollability of node states in a cluster and can be further characterized in an algebraic manner. It is shown in \cite{Ishizaki2014} that the pair $(\Gamma, F)$ can be converted into a positive \textit{tridiagonal realization} by a unitary matrix $T$. Define $\Phi : = - T \Gamma^{-1} F \in \mathbb{R}^{n \times n}$. Then, a cluster $\mathcal{C}_k$ is reducible if and only if there exists a vector $ \phi_k \in \mathbb{R}^{|\mathcal{C}_k|}$ such that
$I({\mathcal{C}_k})   \Phi  = p_k   \phi_k^\top$, where $p_k: = - I({\mathcal{C}_k}) \Gamma^{-1} F $. 
With the vector $p_k$, an aggregation matrix can be defined as  
\begin{equation} \label{eq:aggregation}
\bar{\Pi}: = \begin{bmatrix}
I({\mathcal{C}_1}) & \cdots & I({\mathcal{C}_r})
\end{bmatrix} \blkdiag\left(  \frac{p_1}{\|p_1\|}, ...,\frac{p_r}{\|p_r\|} \right) .
\end{equation}
which can be viewed as a weighted characteristic matrix $\Pi$ in Definition~\ref{defn:cluster}. Note that $\bar{\Pi}^\top \bar{\Pi} = I_r$, and the reduced-order network model then becomes $\hat{g}(s) = \bar{\Pi}  (sI_r - \bar{\Pi}^\top \Gamma \bar{\Pi})^{-1} \bar{\Pi}^\top F$.

If all the clusters are reducible, then the obtained reduced-order network model in Equation~\ref{sysr} has exactly the same input-output behavior as that of the original network, thus $\| g(s) - \hat{g}(s) \|_{\mathcal{H}_2} = 0$. To further reduce the network model, the so-called $\theta$-\textit{reducibility} is defined for a cluster $\mathcal{C}_k$ as  
$\| I({\mathcal{C}_k})   \Phi - p_k   \phi_k^\top \|_\mathrm{F} \leq \theta$, for vectors $p_k, \phi_k \in \mathbb{R}^{|\mathcal{C}_k|}$. If all the clusters are $\theta$-reducible, a \textit{posteriori} upper bound on the reduction error can be formed:
\begin{equation} \label{eq:IshizakiError}
\| g(s) - \hat{g}(s) \|_{\mathcal{H}_2} \leq \gamma   \sqrt{\sum_{k=1}^{r} |\mathcal{C}_k| (|\mathcal{C}_k| - 1) } \theta,
\end{equation}
where $\gamma$ is an upper bound of $\| \bar{\Pi}  (sI_r - \bar{\Pi}^\top \Gamma \bar{\Pi})^{-1} \bar{\Pi}^\top \Gamma + I_n \|_{\mathcal{H}_\infty}$ and characterized by a Riccati inequality.

A generalization of the above method is provided in \cite{ishizaki2015directed} that considers semistable directed networks. In this case, the Frobenius eigenvector of $\Gamma$ is used for constructing the aggregation matrix in Equation~\ref{eq:aggregation} to preserve both semistability and positivity in the reduced-order network model. A so-called projected controllability Gramian, which can be viewed as a special pseudo controllability Gramian, is used for the characterization. The \textit{posteriori} error bound in Equation~\ref{eq:IshizakiError} is also extended to directed networks. However, this extension turns out to be  questionable, since the relevant Riccati inequality in general does not have a solution for semistable systems. A correct formulation in terms of the $\mathcal{H}_\infty$ norm of a semistable system can be found in \cite{ChengAuto2020Gramian}, which has a form of Equation~\ref{eq:LMI}.

The notion of reducibility and the error bound are essential for the clustering procedure in  \cite{Ishizaki2014,ishizaki2015directed}. The core step in the clustering algorithm is to produce a set of $\theta$-reducible clusters, where the value of  $\theta$ is  adjusted in relation to the approximation error bound. This  approach is extended in  \cite{ishizaki2015clustered2o, Ishizaki2016dissipative} to reduce stable second-order network systems, and to reduce networked dissipative systems in the form of Equation~\ref{sysh}.

\subsubsection{$\mathcal{H}_2$-Suboptimal Methods} \label{sec:H2subopt}
Model reduction of a network system can be formulated as nonconvex optimization problem of which the objective function is the $\mathcal{H}_2$ reduction error. The characteristic matrix $\Pi$ is the optimization variable and is subject to the constraint: 
\begin{equation} \label{eq:constraint}
\Pi \in \mathscr{C} : = \left\{\Pi \mathds{1}_r = \mathds{1}_n,~\text{and}~[\Pi]_{ij} = \{0, 1\},\  \forall i = 1,...,n, \; j = 1,..., r.\right\}.
\end{equation}
The optimization problem itself is nonconvex due to the nonlinear objective function in terms of the $\mathcal{H}_2$ norm and the binary variable $\Pi$. In order to solve such a nonconvex problem, a relaxation
of the binary constraints can be taken, leading to suboptimal approaches. 

When we would drop the constraint $\Pi \in \mathscr{C}$ we obtain an $\mathcal{H}_2$ optimal model reduction problem for a generic LTI system which can be solved with the so-called \textit{Iterative Rational Krylov Algorithm} (IRKA)  to seek for a (locally) optimal solution \cite{Petar2015CDC}. The algorithm gives a subspace of dimension $r$ with the basis $V_r \in \mathbb{R}^{n \times r}$. Different from $V$ defined in Equation~\ref{eq:projection}, this $V_r$ is not a feasible solution since it does not preserve the network structure. Therefore, the idea is to find a $\Pi$ matrix in the feasible set $\mathscr{C}$ such that the image of
$\Pi$ is approximately equal to the image of $V_r$, i.e., $\im(\Pi) \approx \im(V_r)$. To this end, \cite{Petar2015CDC} adopts a approach based on a QR decomposition with column pivoting, which is originated from \cite{zha2002spectral} for solving K-means clustering problems.

An alternative method in \cite{ChengECC2020} further studies this nonconvex optimization problem and specifies the objective function using the controllability and observability Gramians. Consider the single-integrator network in Equation~\ref{sys} with $\Gamma = -L$ representing a connected undirected graph. Then, we aim for the following optimization problem:
\begin{align}
\label{eq:optimization_mi}
\min_{\Pi \in {\mathbb R}^{n \times r}} \; J (\Pi): = \tr(B_e^\top P_e B_e)  
\quad \text{s.t.:}  \quad \text{Equation~\ref{eq:constraint}} \ \text{and}\  A_e^\top P_e + P_e A_e + C_e^\top C_e =0,
\end{align}
where $P_e$ is the observability Gramian of the stable system $(A_e, B_e, C_e)$ defined as
\begin{equation} \label{eq_AeBeCe}    
A_e = - 
\begin{bmatrix}
\mathcal{S}_n^\dagger L \mathcal{S}_n & 0 \\ 0 & \mathcal{S}_r^\dagger \Pi^\top L \Pi \mathcal{S}_r
\end{bmatrix}, \  
B_e  =  \begin{bmatrix} \mathcal{S}_n^\dagger F \\ \mathcal{S}_r^\dagger \Pi^\top F \end{bmatrix},
\ C_e = \begin{bmatrix} H \mathcal{S}_n &  - \frac{r}{n} {H} \Pi \mathcal{S}_r\end{bmatrix}.
\nonumber
\end{equation}
The matrices $\mathcal{S}_n$ and $\mathcal{S}_r$ are given by $\mathcal{S}_k=[ I_{k-1} \ \mathds{1}_{k-1}]^\top $ for $k=r,n$. The latter matrices are used to filter out the subspace corresponding to the zero eigenvalues, so that $A_e$ is Hurwitz.  From the objective function $J(\Pi)$, we can derive an explicit expression for its gradient $\nabla J(\Pi)$ so that 
gradient-based algorithms, including projected gradient descent,  Frank-Wolfe optimization, and conditional gradient methods, can be applied to solve the optimization problem in Equation~\ref{eq:optimization_mi}, see \cite{ChengECC2020,lacoste2016FW} and the references therein for more details. 


\subsubsection{Edge Weighting}

Instead of seeking for a way to do the graph clustering, the $\mathcal{H}_2$ optimization scheme can also be applied to construct a ``good'' reduced-order model from a given clustering. To achieve this, we have to go beyond the framework of Petrov-Galerkin projection.   Given a certain clustering, the topology of a reduced-order network is known, while the coupling strengths (edge weights) are considered as free parameters to be determined, \cite{ChengECC2019weight,ChengTAC2020Weight}.

Consider a network system in Equation~\ref{sys} with a connected undirected graph $\mathcal{G}$. Let $\{ \mathcal{C}_1,\mathcal{C}_2,\cdots,\mathcal{C}_r \}$ be a given graph clustering of $\mathcal{G}$. Then, a \textit{quotient graph} $\hat{\mathcal{G}}$ is a $r$-node directed graph obtained by aggregating all the nodes in each cluster as a single node, while retaining connections between clusters and ignoring the edges within each cluster. 
The incidence matrix $\hat{\mathcal{B}}$ of the quotient graph $\hat{\mathcal{G}}$ can be obtained by removing all the zero columns of $\Pi^\top \mathcal{B}$, where $\mathcal{B}$ is the incidence matrix of $\mathcal G$, and $\Pi$ is the characteristic matrix of the clustering.
Denote 
$
\hat{W} = \diag(\hat{w})$,  {with} 
$\hat{w} = [
\hat{w}_1 \ \hat{w}_2 \ ... \ \hat{w}_\kappa 
]^\top, 
$
as the edge weighting matrix of $\hat{\mathcal{G}}$, with $\hat{w}_k \in \mathbb{R}_{>0}$ and $\kappa$ the number of edges in $\hat{\mathcal{G}}$. Then, the parameterized model of a reduced-order network is obtained as
\begin{equation} \label{sysrp}
(\Pi^\top \Pi)\dot{z}(t) = - \hat{\mathcal{B}} \hat{W} \hat{\mathcal{B}}^\top z(t) +  \Pi^\top {F} u (t), \quad
\hat{y}(t) = H \Pi z(t),
\end{equation} 
which has as transfer matrix $\hat{\eta}_p (s, \hat{W}) = H \Pi (s \Pi^\top \Pi + \hat{\mathcal{B}} \hat{W} \hat{\mathcal{B}}^\top)^{-1} \Pi^{-1} F$. In the reduced-order model, the edge weight matrix $\hat{W}$ is the only unknown, and can be determined via an optimization procedure. 

\begin{figure}[t]
	\centering
	\includegraphics[width=4.2in]{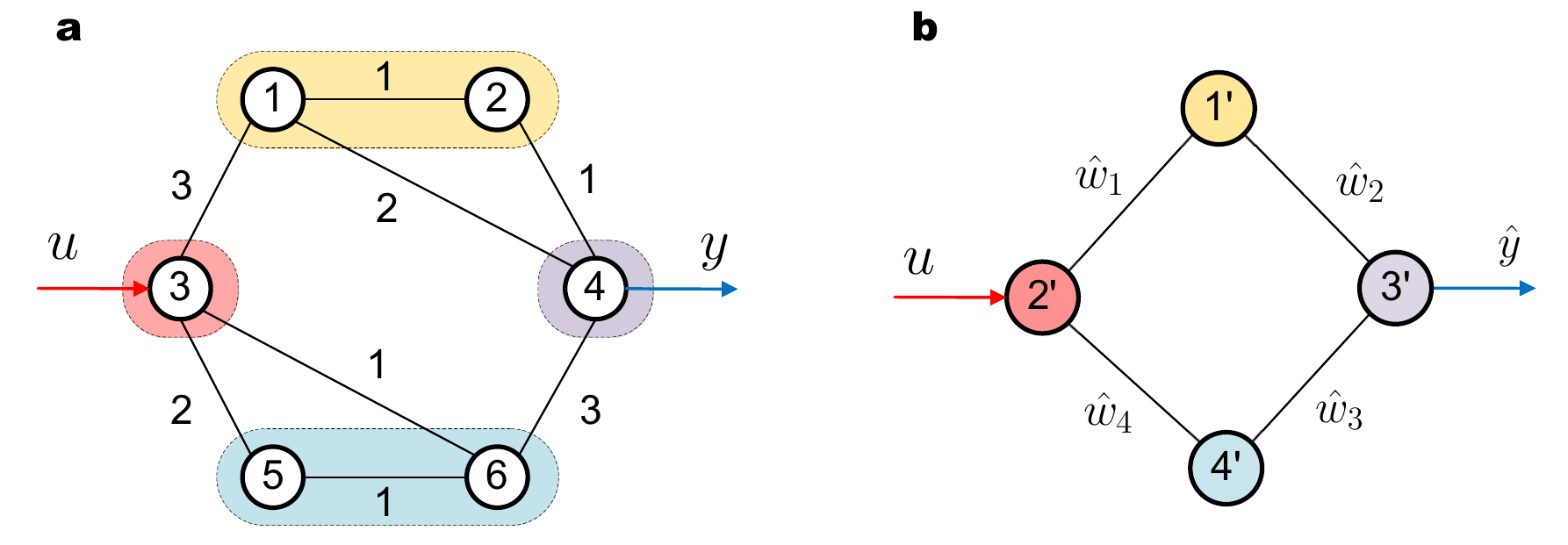}   
	\caption{(\textit{a}) An undirected network consisting of 6 nodes, in which node 3 is the leader and node 4 is measured. Four clusters are indicated by different colors. (\textit{b}) A quotient graph obtained by clustering.}
	\label{fig:edgeweighting}
\end{figure}

\begin{exm}
	Consider an undirected graph composed of 6 nodes in \textbf{Figure~\ref{fig:edgeweighting}a}. Given a clustering with 
	$\mathcal{C}_1 = \{1,2\},~ \mathcal{C}_2 = \{3\},~ \mathcal{C}_3 = \{4\},~ \mathcal{C}_4 = \{5,6\},$ the quotient graph is obtained in \textbf{Figure~\ref{fig:edgeweighting}b} with 
	edge weight matrix $\hat{W} = \diag(\hat{w}_1, \hat{w}_2, \hat{w}_3, \hat{w}_4)$. Then, the parameterized model of the reduced network is constructed as
	\begin{align*} 
	\underbrace{\begin{bmatrix}
		2 & 0 & 0 & 0\\
		0 & 1 & 0 & 0\\
		0  & 0 & 1 & 0\\
		0  & 0 & 0 & 2
		\end{bmatrix}}_{\Pi^\top \Pi}
	\begin{bmatrix}
	\dot{z}_1 \\ \dot{z}_2 \\\dot{z}_3 \\ \dot{z}_4
	\end{bmatrix}
	= 
	-\underbrace{\begin{bmatrix}
		\hat{w}_1 + \hat{w}_2 & - \hat{w}_1 & - \hat{w}_2 & 0\\
		-\hat{w}_1 & \hat{w}_1 + \hat{w}_4 & 0 & - \hat{w}_4\\
		- \hat{w}_2 & 0 & \hat{w}_2 + \hat{w}_3 & -\hat{w}_3\\
		0 & - \hat{w}_4 & - \hat{w}_3 & \hat{w}_3 + \hat{w}_4\\
		\end{bmatrix}}_{	\hat{\mathcal{B}} \hat{W} \hat{\mathcal{B}}^\top} 
	\begin{bmatrix}
	{z}_1 \\ {z}_2 \\{z}_3 \\ {z}_4 
	\end{bmatrix}
	+
	\underbrace{\begin{bmatrix}
		0  \\
		1 \\
		0 \\
		0
		\end{bmatrix}}_{\Pi^\top F} u,  
	\end{align*}
	and $ 
	y 
	= {z}_3$. In this model, the diagonal elements of $\hat{W}$ are the parameters to be determined.
\end{exm}

Similar to the process in Section~\ref{sec:H2subopt}, an optimization problem can be formulated  to minimize the reduction error between the original and reduced-order network systems by tuning the edge weights. Specifically, the objective function is $\| \eta(s) - \hat{\eta}_p(s, \hat{W}) \|_{\mathcal{H}_2}$, in which a diagonal and positive definite $\hat{W}$ is the optimization variable.

There are multiple algorithms for solving such a problem. In \cite{ChengECC2019weight}, the upper bound of the reduction error, i.e., the expression $\| \eta(s) - \hat{\eta}_p(s, \hat{W}) \|_{\mathcal{H}_2} < \gamma$, is characterized by a set of LMIs, and a cross iteration algorithm is applied such that the upper bound $\gamma$ decreases via iterations. An alternative approach is presented in \cite{ChengTAC2020Weight} by means of a \textit{convex-concave decomposition}. This approach is inspired by the work in \cite{dinh2011combining}, and is based on linearization of the optimization problem at a given point, so that the problem becomes \textit{convex} and can be solved efficiently.  Then, the overall problem can be solved in an iterative fashion, and in each iteration a convex optimization problem needs to be solved.

It is worth noting that the edge weighting approach can be implemented as a subsequent procedure for the  clustering-based methods described earlier, i.e., 
we can first apply an aforementioned algorithm to find a graph clustering, whose edge weights can be used to initialize $\hat{W}$ in the edge weighting approach. Then, through iterations, a more accurate reduced-order network model can be generated.

\subsection{Other Topological Reduction Methods}

There exists a vast amount of literature
about the problem of topological reduction. 
In this section, we briefly summarize several other representative methods. 

\subsubsection{Singular Perturbation Approximation and Kron reduction}
Along with graph clustering, the other mainstream methodology for simplifying the topological complexity of a network is based on time-scale separation analysis, and in particular, singular perturbation approximation \cite{kokotovic1976singular}. This method has been extensively investigated in the applications of biochemical reaction systems \cite{Rao2013graph,anderson2011reduction,hancock2015simplified} and electric power networks \cite{chow1995inertial,Biyik2008Areaaggregation,Chow2013PowerReduction,romeres2013novel}. A key feature of those systems is that there usually is an explicit or non-explicit separation of time scales in the states of networks. For example, slow coherency theory implies that power networks are naturally decomposed into areas, where power generators within each area synchronize on a faster time scale, while network-wide interactions  between the areas are captured by slower motions. Singular perturbation methods then help to separate these dynamics to analyze them separately.  This produces a reduced-order network model that retains for example the low frequency behavior of a large-scale network. How to identify and separate fast/slow states is a crucial step in this method, and its application is highly dependent on the specific application.

A relevant concept is \textit{Kron reduction} of graphs, which is a terminology commonly used in classic circuit theory and in related fields such as electrical impedance tomography, and transient stability assessment in power networks, see e.g., \cite{Chow2013PowerReduction,Dorfler2013Kron,Monshizadeh2017Constant}. This may also be used for exact reduction, e.g., to go from a differential algebraic description of the a network system to a differential description {with structure preservation}.

%


\subsubsection{Kullback-Leibler Aggregation} 

Networked dynamical systems derived from the discretization of thermodynamics and fluid dynamics are usually modeled as regular discrete-time Markov chains without control inputs. For this type of systems, the notion of \textit{Kullback-Leibler (K-L) divergence rate} can be adopted to measure the difference between two Markov chains, defined on the same state space. {This notion is further extended in \cite{deng2011optimal} to measure the K–L divergence rate between the original and reduced-order models defined on different state spaces. 
With the new K–L divergence rate,} an optimization problem is formulated which aims to find an optimal partition of the states, which are aggregated to form a reduced-order model.  
An application of this method is explored in reducing complex building thermal systems \cite{deng2014structure}. 

\subsubsection{Network Reduction towards Scale-Free Structure}

Graph clustering and aggregation have also been studied to retain more relevant properties such as connectivity and scale-freeness in \cite{Martin2018scalefree,Martin2018mergetocure}. The scale-free networks typically contain a few nodes with a large degree (the so-called hubs) and a large number of nodes with a small degree, and the distribution of the node degree follows a certain power law \cite{Martin2018scalefree}.

The model reduction problem of networks preserving the scale-free property can be formulated as an optimization problem that finds a clustering of a given large-scale network such that the aggregated network has a degree distribution closest to a desired scale-
free distribution.  The preservation of the scale-free structure is particularly important for applications of flow networks, such as traffic networks, power networks or packet flow networks.

\subsubsection{Indirect Network Reduction Methods}
\label{sec:indirect}
Different from the mechanisms of clustering and aggregation that  directly produce a reduced network, indirected methods in \cite{ChengIFAC2017,ChengAuto2019,LanlinCDC2019second-order,LanlinCDC2019eigenvalue} seek a structure-preserving reduced-order model using a two-step procedure: reduction and transformation. In the reduction step, a lower-dimensional model of a given large-scale network system is constructed by using conventional model reduction methods, e.g. generalized balanced
truncation \cite{ChengIFAC2017,ChengAuto2019} or moment matching \cite{LanlinCDC2019second-order}. Generally, the reduced-order model generated in this step does not allow for a network interpretation. Then, the transformation step is implemented subsequently, which converts the reduced-order model obtained in the previous step into a network model. This method is also relevant for the combined nodal and topological reduction procedure in Section \ref{sec:combine}.

A key for the transformation is presented in  \cite{ChengIFAC2017,ChengAuto2019}: A matrix is similar to a Laplacian matrix of a connected undirected graph if and only if it is diagonalizable and has exactly one zero eigenvalue while all the other eigenvalues are real positive. This result guides the reduction step, in which certain spectral constraint has to be imposed. Then, the second step turns out to be an eigenvalue matching problem.

Similarly, an eigenvalue assignment approach,
\cite{LanlinCDC2019eigenvalue}, directly selects a subset of the Laplacian spectrum of the original network to be the eigenvalues of the Laplacian matrix for the reduced network. By doing so, certain properties of original network such as stability and synchronization can be preserved through the reduction process.


\section{REDUCTION OF THE FULL SYSTEM DYNAMICS} 

In the previous section the reduction of the topological
structure based on clustering methods is considered. The input-output 
structure is considered when looking at for example the ${\cal H}_2$ norm, but 
for control systems it is undoubtedly very important that the model preserves  
certain input-output or control structures beyond only considering the ${\cal H}_2$ norm. For general linear and 
nonlinear systems various methods are developed, e.g., 
\cite{Antoulas05}, and \cite{Fujimoto2010BT}, \cite{Scherpen11}. 
We refer to the Sidebar on Linear Systems for a very brief introduction. Here we first 
consider reduction of the nodal dynamics while preserving certain
graph properties. Secondly, we treat a combined nodal and 
topological method.  

\subsection{Reduction of the Nodal Dynamics}

\label{sec:nodalredution}

For the nodal dynamics it is useful to consider reduction methods for interconnected systems. In the corresponding literature the network perspective is not the primary focus, but the methods are nevertheless relevant for network systems and treated in this subsection.  In addition, we treat reduction methods for nodal dynamics that do explicitly take the network perspective and preserve properties like synchronization.  

\subsubsection{Reduction Methods for Interconnected Linear Systems}
\begin{figure}[t]
	\centering
	\includegraphics[width=6cm]{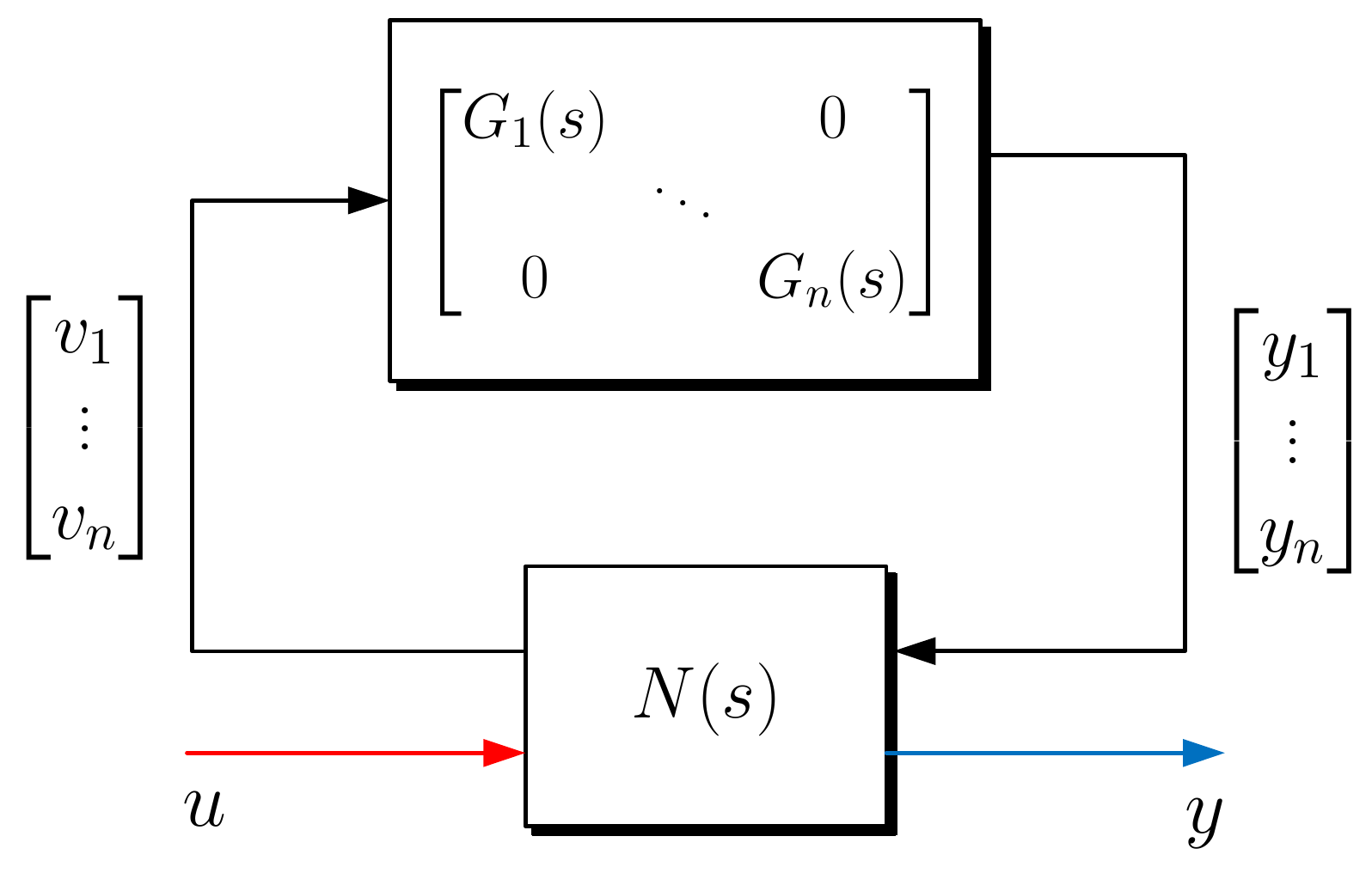}
	\caption{The interconnected system with subsystems to be reduced, $G_1(s), \ldots, G_n(s)$ stored in a block-diagonal transfer matrix $G(s)$ and connected to the transfer matrix $N(s)$ that models the interconnection topology and excitation and measurement dynamics, \cite{sandberg2009interconnected}.}
	\label{fig:sandberg}
\end{figure}

Perhaps one of the first relevant papers that considers reduction methods for interconnected systems is  \cite{sandberg2009interconnected}, where linear fractional transformations are considered, see \textbf{Figure~\ref{fig:sandberg}}. 
The latter paper considers two methods for reduction, both based on balanced truncation principles. The first method only considers the diagonal blocks of the observability and controllability Gramians, with each block corresponding to a subsystem, hence neglecting the off diagonal blocks. An expression for the posteriori error bound is provided, but unfortunately, an a priori error bound is not obtained. The second method uses generalized Gramians, i.e., instead of considering the observability and controllability Gramians that are the unique solutions to corresponding Lyapunov equations, the (non-unique) solutions to  Lyapunov inequalities are considered. The freedom in choosing solutions to these Lyapunov inequalities provides a possibility to pick block diagonal solutions, and consequently results in an a priori error bound. An extension towards a singular perturbation perspective is provided in \cite{ishizaki2013structured}, and balanced truncation based on generalized Gramians in the discrete time setting for interconnected systems is provided in \cite{jaoude2017balanced}. More generally, reduction methods for coupled systems, and an overview of the various methods until 2008 are provided in \cite{Vandendorpe2008} and  \cite{Reis2008survey}.

A more recent result deals with interconnected systems in a graph setting, e.g., \cite{Ionescu2020H2}. In particular, the subsystems are of the form
$\dot{x}_i = \sum_{j \in {\cal N}_i} A_{ij} x_j + B_i u, 
$ 
where ${\cal N}_i$ is the index set of the connections of the $i^{th}$ subsystem with other subsystems. With help of generalized Gramians, optimal ${\cal H}_2$ moment matching based model reduction is performed while preserving the network structure. This is done in a convex optimization setting, and additionally a projected gradient method is applied for the non-convex case.

\subsubsection{Nodal Reduction while Preserving Synchronization Properties}
The above methods consider interconnected systems, but do not take into account yet the preservation of properties that are typically relevant for (controlled) network systems, such as consensus or synchronization properties. For that, additional steps have to be taken. To the best of our knowledge the first work that considers preservation of synchronization properties in network systems is \cite{Monshizadeh2013stability}. Consider a diffusively-coupled linear network system as in Equation~\ref{sysh} with a symmetric Laplacian $\Gamma = -L$. Then a bounded real balancing method can be employed in order to preserve the stability of the network. In addition, synchronization can be preserved by considering the positive definite solution $K$ of the following Riccati equation:
\[
(A- \lambda BC)^T K + K (A -\lambda BC) + C^TC +\delta^2 K BB^TK=0
\]
where $\delta$ is a scalar that has to fulfill some additional conditions, and $\lambda$ is an eigenvalue of the Laplacian $L$ for which $A-\lambda BC$ is Hurwitz. The maximum and minimum solutions
of the equation can now be balanced, i.e., simultaneously diagonalized, similar as is done for standard balancing. Truncating the system based on these diagonal values results in reduced order dynamics of the agents and a network system which is still synchronized. Furthermore, an a priori error bound is provided based on the truncated diagonal values, $\delta$, and the largest eigenvalue of the Laplacian.

The above balancing method based on the minimum and maximum solution of a Riccati equation can be generalized to finding a solution to the inequality, and as such this can be related to solutions of Linear Matrix Inequalities (LMIs). Also, so far, only linear dynamics in the nodes is considered, where in practice nonlinearities play an important role. Because error dynamics are more difficult to handle in the case of nonlinear systems, it useful to consider nodal dynamics represented by Lur'e systems, i.e., systems with a static nonlinearity in the feedback loop. In \cite{ChengECC2018Lure,ChengEJC2018Lure} robust synchronization preserving reduction methods for nodal Lur'e systems are considered. For this, consider the following nodal dynamics:
\begin{equation}  
\label{eq:Luresystem}
\bm{\Sigma_i} :  
\left\{
\begin{array}{l}
\dot{x}_i(t) = A x_i(t) + B u_i(t) +  E z_i(t) \\
y_i(t) = C x_i(t) \\
z_i(t) = -\phi (y_i(t))
\end{array}
\right.
\end{equation}
where $x_i \in \mathbb{R}^n$, $u_i,y_i,z_i \in \mathbb{R}^m$, and $\phi:\mathbb{R}^m \rightarrow \mathbb{R}^m$ fulfilling some sector bound condition. Without loss of generality, we take $m=1$, \cite{ChengECC2018Lure}. For a diffusively coupled network which is robustly synchornized and which has nodal dynamics as in Equation \ref{eq:Luresystem}, we now take the minimum and maximum solutions $\mathcal{K}_{min}$ and $\mathcal{K}_{max}$ of the following LMI
\begin{equation*} 
\begin{bmatrix}
{ A}^T \mathcal{K} + \mathcal{K} {A} +  {C^TC} & { \lambda_n} \mathcal{K} { B} & \mathcal{K} {E} - {C}^T \\
{ \lambda_n} { B}^T \mathcal{K} & -I & 0 \\
{E}^T \mathcal{K} - {C} & 0 & 0
\end{bmatrix} \preccurlyeq 0,
\end{equation*}
with $\lambda_n$ the largest eigenvalue of the Laplacian. $\mathcal{K}_{min}$ and $\mathcal{K}_{max}$ can be balanced and reduction based on them results in a robustly synchronized network of Lur'e systems with a priori determined error bounds. Variations on this can be done to obtain even better error bounds, and an extension to the Multi-Input Multi-Output case is provided in \cite{ChengEJC2018Lure}.




\begin{textbox}[h]\section{Model Reduction for Linear Systems}
	It is generally accepted that model reduction approaches for linear control systems can be roughly divided into two types of approaches, i.e., singular value based and moment matching based methods, e.g., \cite{Antoulas05}. A very brief review follows.  Let $(A,B,C)$ be a state space realization of $\Sigma$ with dimension $n$,
	with input $u \in \mR^m$, state $x \in \mR^n$ and output $y \in \mR^p$.
	We assume that the system is {\em asymptotically stable and
		minimal}, i.e., controllable and observable. The corresponding transfer matrix is given by $G(s) = C(sI - A)^{-1}B$.
	
	\subsection{Balanced Realizations}
	\begin{thm} \cite{Moore81} \label{tmtuslin1} 
		Consider the system $\Sigma$. Take $P$ the controllability Gramian and $Q$ the observability Gramian. The eigenvalues of $QP$ are similarity invariants, i.e., they do
		not depend on the choice of the
		sate space coordinates. Furthermore, there exists a state space representation where
		$	\Sigma := Q = P = \diag{\left\{ \sigma_1, \ldots,\sigma_n\right\} }$, 
		with  $\sigma_1 \geq \sigma_2 \geq ... \geq \sigma_n >0$ the square roots of the eigenvalues of $QP$.
		Such representations are called
		{\em balanced}, and the system
		is in {\em balanced form}.  Furthermore, the $\sigma_i$'s, i=1,..,n,
		equal the Hankel singular values,
		i.e., the singular values of the Hankel operator. 
	\end{thm}

	\noindent 
	The Hankel singular values form a measure for the contribution to minimality of a state component. This provides the basis for model reduction methods based on balanced realizations. Model reduction based on balancing is possible with  a-priori error bounds in various norms, e.g.,the ${\mathcal H}_\infty$, and Hankel norm. In particular, it can be shown that balanced truncation results in an error bound corresponding to the sum of the truncated Hankel singular values, i.e., $\|G(s) - G_r(s)\|_{\infty}\leq 2 \sum_{k=r}^{n} \sigma_k$, where $G_r(s)$ represents the transfer matrix of the reduced order system, \cite{glover1984all}.  We refer to \cite{Antoulas05} and \cite{Scherpen11} for a more elaborate overview.
	
	\subsection{Moment Matching}
	The principle of moment matching for a linear system is based on the series representation of the transfer matrix of the system. For more detailed expositions, we refer to e.g., \cite{Antoulas05} and \cite{Scarciotti17}. Without loss of generality it is assumed that $m=p=1$.  
	\begin{defn}\label{defmoment}
		The 0-moment at $s_1 \in \mathbb{C}$ of $\Sigma$ is the complex number $\eta_0(s_1) = C(s_1I - A)^{-1}B$. The k-moment at $s_1$ is given by the complex number 
		\[
		\eta_k(s_1) = \frac{(-1)^k}{k!} \frac{d^k \left\{C(sI-A)^{-1}B\right\}}{ds^k}_{\displaystyle \big|s=s_1} , \,\,\, k=1,2,3,\ldots
		\]
	\end{defn}
	The point $s_1 \in \mathbb{C}$ is called an interpolation point. The
	approximation problem for system $(A,B,C)$ at 
	$s_1 \in \mathbb{C}$ is to find a system $(\hat{A},\hat{B} ,\hat{C})$ of order $\nu<n$,  with transfer function $\hat{G}(s)= \hat{C}(sI - \hat{A})^{-1} \hat{B}$, such that $\eta_k(s_1) = \hat{\eta}_k(s_1)$,  with $\hat{\eta}_k(s_1)$ the moments of $\hat{G}(s)$,
	$k = 1, \ldots, \nu$. Various types of moment matching methods are developed. Generally it is {\it not} possible to provide an a-priori error bound. These methods are however computationally very interesting if one handles systems with millions of states.
	
\end{textbox}

\subsection{Combined Topological and Nodal Reduction}
\label{sec:combine}

\begin{figure}[h]\centering
	\centering
	\includegraphics[width=0.9\textwidth]{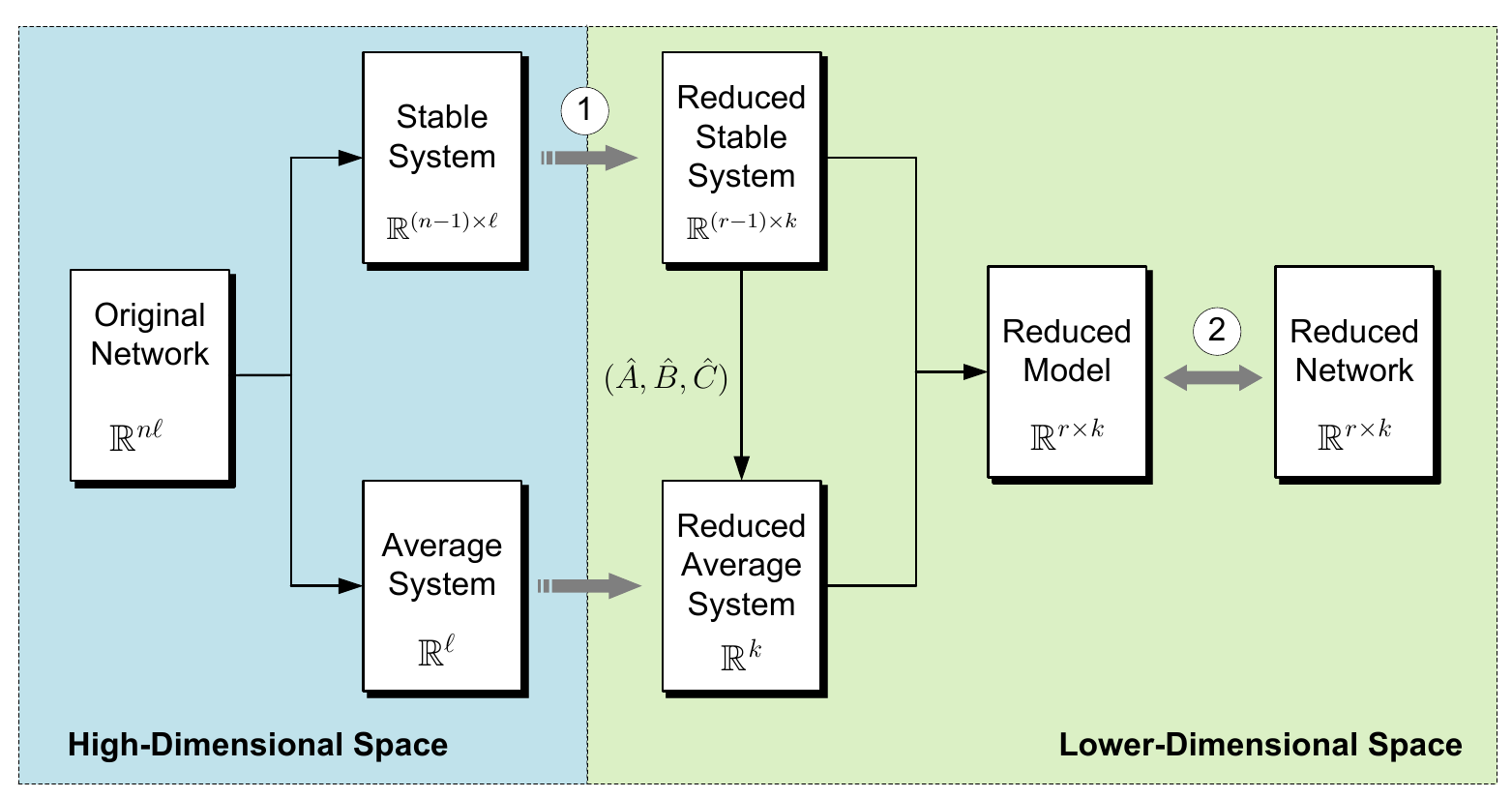}	
	\caption{The model reduction procedure for a networked passive system, which contains two key techniques: balanced truncation based on generalized Gramians and a coordinate transformation that reconstructs a network.}	
	\label{fig:procedure}
\end{figure}

From the previous sections, it has been observed that the techniques for topological simplification and subsystem reduction in network systems are derived from rather different perspective.
The methods for reducing subsystem dynamics are commonly adapted from classic model reduction techniques e.g., balanced truncation, or Krylov subsystem methods, while in structure-preserving simplification of network structures, clustering-based approaches have demonstrated a superior performance. 
In this section, we discuss the combination of the two approximation problems in a unified framework. This is particularly needed when dealing with a network with both complex topology and high-order subsystems.

However, how to perform a simultaneous reduction of topological complexity and subsystem dynamics is not straightforward.
Naively, we may apply the methods of Section~\ref{sec:topology} and Section~\ref{sec:nodalredution} separately to achieve the two approximation goals one by one. Nevertheless, which reduction sequence gives a better approximation is still unclear, and moreover there is hardly a guarantee on the approximation error.  
Relevant results in the literature on combined topological and nodal reduction are rare. The existing ones are developed generally for networked homogeneous linear systems as in Equation~\ref{sysh} under specific assumptions \cite{Ishizaki2016dissipative,ChengAuto2019}. For example, regularity and dissipativity of the entire system matrix $I \otimes A + \Gamma \otimes BC$ is imposed in \cite{Ishizaki2016dissipative}, which admits a block-diagonal Lyapunov function. It is essential for preserving the stability of the reduced-order model and to derive an upper bound on the approximation error caused by the reduction of both network structure and subsystems, where the topological reduction is done by graph clustering.

In \cite{ChengAuto2019}, a network system is considered which is synchronized but not necessarily stable. The synchronization property is based on the assumption of each subsystem  $(A, B, C)$ being minimal and passive. 
Extending the results for networked single integrators of Section~\ref{sec:indirect}, the topological complexity is reduced in an indirect manner. A generalized balanced truncation method then results in a unified framework to simplify the network structure and subsystem
dynamics simultaneously. The reduction scheme is illustrated in \textbf{Figure~\ref{fig:procedure}}, which contains two core steps.

The first step is to decompose the considered network system into two parts which correspond to the average of all subsystems and the discrepancy among the subsystem states, respectively. Due to the synchronization property, the latter part that captures the main dynamics of the entire network is represented by an asymptotically stable system $\bf{\Sigma_s}$ with a Hurwitz system matrix $(I_{n-1} \otimes A - \bar{\Lambda} \otimes BC)$, where $\bar{\Lambda} \in \mathbb{R}^{(n-1) \times (n-1)}$ shares all the nonzero eigenvalues with the Laplacian matrix $L \in \mathbb{R}^{n \times n}$ of the original network. Then, a pair of generalized Gramians with a Kronecker product structure are selected, which is key to decouple the balanced truncation procedures of the network structure and subsystem dynamics. In the second step, the simplified stable system and average system are integrated, resulting in a reduced-order model of the original $\bm{\Sigma}$. However, such reduced-order model only provides an approximation of the input-output mapping of $\bm{\Sigma}$, rather than the network structure. To restore the network interpretation in that reduced-order model, a coordinate transformation is required to recover a reduced Laplacian matrix. A theoretical foundation for such transformation is firstly presented in \cite{ChengIFAC2017}:
\begin{thm} \label{thm:LapReal}
	A real square matrix $\mathcal{N}$ is similar to the Laplacian matrix ${L}$ associated with a weighted undirected connected graph, if and only if $\mathcal N$ is diagonalizable and has an eigenvalue at $0$ with multiplicity $1$ while all the other eigenvalues are real and positive.
\end{thm}
In \cite{ChengAuto2019}, a detailed proof of the above theorem is provided, and meanwhile a method is given for reconstructing an undirected network from  given eigenvalues $0 = \lambda_1 < \lambda_2 \leq \cdots \leq \lambda_{n}$. Recently, an alternative graph reconstruction approach is presented in \cite{forrow2018spectra}. The reconstructed graphs are usually complete, i.e., there is an edge between any pair of nodes. A subsequent procedure can be taken to sparsify the interconnection structure, see e.g., \cite{forrow2018spectra,jongsma2017CycleRemoval}.



\begin{summary}[SUMMARY POINTS]
	\begin{enumerate}
		\item Graph-theoretical analysis plays a paramount role in reduced-order modeling of complex network systems.
		
		\item Clustering methods provide a structure preserving manner to reduce the topology of network systems. Dissimilarity-based clustering provides a rather general framework for simplifying the topological complexity of a network system, where the key is to properly define a metric to characterize the dissimilarity between nodes/clusters.
		
		\item Network systems that achieve synchronization are naturally semistable. Novel pseudo Gramian notions are introduced for semistable systems as an extension of Gramian matrices for asymptotically stable systems, providing a useful tool to characterize dissimilarity and compute the reduction error. 
		
		\item Generalized balanced truncation allows more freedom in constructing a reduced-order model with some desired structures/properties. Thus, it is widely used in the approximation of network systems, particularly in dealing with subsystem reduction.
		
		\item Graph reconstruction can realize a network representation from a reduced-order model satisfying certain spectral constraints. This makes it possible to apply more classical model order reduction methods from the control systems literature for structure preserving reduction of network systems.  
	\end{enumerate}
\end{summary}

\begin{issues}[FUTURE ISSUES]
	\begin{enumerate}
		\item  The approximation of complex network systems with nonlinear couplings and nonlinear subsystems is still a challenge, requiring further investigation.
		\item How to reduce the topological complexity of dynamic networks composed of heterogeneous subsystems is not clear yet.
		\item The application of reduced-order network models for designing controllers and observers for large-scale networks is appealing, while obtaining provable guarantees on the functionality of the controllers/observers based on reduced-order models should be further
		explored.
	\end{enumerate}
\end{issues}


\section*{ACKNOWLEDGMENTS}
Posted with permission from the Annual Review of Control, Robotics, and Autonomous Systems, Volume 4 by Annual Reviews, http://www.annualreviews.org.

\bibliographystyle{ar-style3} 
\bibliography{ref} 

\begin{thebibliography}{102}
\expandafter\ifx\csname natexlab\endcsname\relax\def\natexlab#1{#1}\fi

\bibitem{mesbahi2010graph}
Mesbahi M, Egerstedt M. 2010.
\textit{Graph Theoretic Methods in Multiagent Networks}.
Princeton University Press

\bibitem{Newman2003Review}
Newman ME. 2003.
The structure and function of complex networks.
\textit{SIAM Review} 45:167--256

\bibitem{ren2005survey}
Ren W, Beard RW, Atkins EM. 2005.
\textit{A survey of consensus problems in multi-agent coordination}.
In \textit{Proceedings of the 2005 American Control Conference (ACC)}, pp.
  1859--1864. IEEE

\bibitem{Rao2013graph}
Rao S, van~der Schaft AJ, Jayawardhana B. 2013.
A graph-theoretical approach for the analysis and model reduction of
  complex-balanced chemical reaction networks.
\textit{Journal of Mathematical Chemistry} 51:2401--2422

\bibitem{Ahsendorf2014GeneRegulation}
Ahsendorf T, Wong F, Eils R, Gunawardena J. 2014.
A framework for modelling gene regulation which accommodates non-equilibrium
  mechanisms.
\textit{BMC Biology} 12:102

\bibitem{proskurnikov2017tutorial}
Proskurnikov AV, Tempo R. 2017.
A tutorial on modeling and analysis of dynamic social networks. part i.
\textit{Annual Reviews in Control} 43:65--79

\bibitem{dorfler2012synchronization}
Dorfler F, Bullo F. 2012.
Synchronization and transient stability in power networks and nonuniform
  kuramoto oscillators.
\textit{SIAM Journal on Control and Optimization} 50:1616--1642

\bibitem{Chow2013PowerReduction}
Chow JH. 2013.
\textit{Power System Coherency and Model Reduction}.
Springer

\bibitem{Antoulas05}
Antoulas AC. 2005.
\textit{Approximation of Large-Scale Dynamical Systems}.
Philadelphia, USA: SIAM

\bibitem{Bai2002Krylov}
Bai Z. 2002.
Krylov subspace techniques for reduced-order modeling of large-scale dynamical
  systems.
\textit{Applied Numerical Mathematics} 43:9--44

\bibitem{Astolfi2010Moment}
Astolfi A. 2010.
Model reduction by moment matching for linear and nonlinear systems.
\textit{IEEE Transactions on Automatic Control} 55:2321--2336

\bibitem{Moore81}
Moore BC. 1981.
Principal component analysis in linear systems: Controllability, observability,
  and model reduction.
\textit{IEEE Transactions on Automatic Control} 26:17--32

\bibitem{glover1984all}
Glover K. 1984.
All optimal {Hankel}-norm approximations of linear multi-variable systems and
  their ${L}_\infty$-error bounds.
\textit{International Journal of Control} 39:1115--1193

\bibitem{besselink2013comparison}
Besselink B, Tabak U, Lutowska A, Van~de Wouw N, Nijmeijer H, et~al. 2013.
A comparison of model reduction techniques from structural dynamics, numerical
  mathematics and systems and control.
\textit{Journal of Sound and Vibration} 332:4403--4422

\bibitem{summers2014optimal}
Summers TH, Lygeros J. 2014.
Optimal sensor and actuator placement in complex dynamical networks.
\textit{IFAC Proceedings Volumes} 47:3784--3789

\bibitem{gates2016control}
Gates AJ, Rocha LM. 2016.
Control of complex networks requires both structure and dynamics.
\textit{Scientific Reports} 6:1--11

\bibitem{Kim2018role}
Kim JZ, Soffer JM, Kahn AE, Vettel JM, Pasqualetti F, Bassett DS. 2018.
Role of graph architecture in controlling dynamical networks with applications
  to neural systems.
\textit{Nature Physics} 14:91--98

\bibitem{ishizaki2018graph}
Ishizaki T, Chakrabortty A, Imura JI. 2018.
Graph-theoretic analysis of power systems.
\textit{Proceedings of the IEEE} 106:931--952

\bibitem{Obinata2012MR}
Obinata G, Anderson BDO. 2012.
\textit{Model reduction for Control System Design}.
Springer Science \& Business Media

\bibitem{Mustafa1991Controller}
{Mustafa} D, {Glover} K. 1991.
Controller reduction by ${H}_\infty$ balanced truncation.
\textit{IEEE Transactions on Automatic Control} 36:668--682

\bibitem{SurveyClustering}
Jain AK, Murty MN, Flynn PJ. 1999.
Data clustering: a review.
\textit{ACM Computing Surveys (CSUR)} 31:264--323

\bibitem{Schaeffer2007SurveyClustering}
Schaeffer SE. 2007.
Graph clustering.
\textit{Computer Science Review} 1:27--64

\bibitem{Schaft2014}
van~der Schaft AJ. 2014.
\textit{{On model reduction of physical network systems}}.
In \textit{Proceedings of 21st International Symposium on Mathematical Theory
  of Networks and Systems (MTNS)}, pp.  1419--1425. Groningen, The Netherlands

\bibitem{Ishizaki2014}
Ishizaki T, Kashima K, Imura JI, Aihara K. 2014.
{Model reduction and clusterization of large-scale bidirectional networks}.
\textit{IEEE Transactions on Automatic Control} 59:48--63

\bibitem{Monshizadeh2014}
Monshizadeh N, Trentelman HL, Camlibel MK. 2014.
Projection-based model reduction of multi-agent systems using graph partitions.
\textit{IEEE Transactions on Control of Network Systems} 1:145--154

\bibitem{deng2014structure}
Deng K, Goyal S, Barooah P, Mehta PG. 2014.
Structure-preserving model reduction of nonlinear building thermal models.
\textit{Automatica} 50:1188--1195

\bibitem{Besselink2016Clustering}
Besselink B, Sandberg H, Johansson KH. 2016.
Clustering-based model reduction of networked passive systems.
\textit{IEEE Transactions on Automatic Control} 61:2958--2973

\bibitem{ChengTAC20172OROM}
Cheng X, Kawano Y, Scherpen JMA. 2017.
Reduction of second-order network systems with structure preservation.
\textit{IEEE Transactions on Automatic Control} 62:5026 -- 5038

\bibitem{ChengTAC2018MAS}
Cheng X, Kawano Y, Scherpen JMA. 2019.
Model reduction of multi-agent systems using dissimilarity-based clustering.
\textit{IEEE Transactions on Automatic Control} 64:1663--1670

\bibitem{jongsma2018model}
Jongsma HJ, Mlinari{\'c} P, Grundel S, Benner P, Trentelman HL. 2018.
Model reduction of linear multi-agent systems by clustering with $h_2$ and
  ${H}_\infty$ error bounds.
\textit{Mathematics of Control, Signals, and Systems} 30:6

\bibitem{Martin2018scalefree}
Martin N, Frasca P, Canudas-de Wit C. 2018.
Large-scale network reduction towards scale-free structure.
\textit{IEEE Transactions on Network Science and Engineering} 6:711--723

\bibitem{Reis2008survey}
Reis T, Stykel T. 2008.
A survey on model reduction of coupled systems.
In \textit{Model order reduction: theory, research aspects and applications}.
  Springer

\bibitem{Vandendorpe2008}
Vandendorpe A, Van~Dooren P. 2008.
Model reduction of interconnected systems.
In \textit{Model order reduction: theory, research aspects and applications}.
  Springer

\bibitem{sandberg2009interconnected}
Sandberg H, Murray RM. 2009.
Model reduction of interconnected linear systems.
\textit{Optimal Control Applications and Methods} 30:225--245

\bibitem{Monshizadeh2013stability}
Monshizadeh N, Trentelman HL, Camlibel MK. 2013.
Stability and synchronization preserving model reduction of multi-agent
  systems.
\textit{Systems \& Control Letters} 62:1--10

\bibitem{ChengEJC2018Lure}
Cheng X, Scherpen JMA, Zhang F. 2019.
Model reduction of synchronized homogeneous {Lur'e} networks with incrementally
  sector-bounded nonlinearities.
\textit{European Journal of Control} 50:11--19

\bibitem{Ishizaki2016dissipative}
Ishizaki T, Ku R, Imura Ji. 2016.
\textit{Clustered model reduction of networked dissipative systems}.
In \textit{Proceedings of the 2016 American Control Conference}, pp.
  3662--3667. IEEE

\bibitem{ChengAuto2019}
Cheng X, Scherpen JMA, Besselink B. 2019.
Balanced truncation of networked linear passive systems.
\textit{Automatica} 104:17--25

\bibitem{Bullo2019Lectures}
Bullo F. 2019.
\textit{Lectures on Network Systems}.
Kindle Direct Publishing

\bibitem{Wu2007Synchronization}
Wu CW. 2007.
\textit{Synchronization in Complex Networks of Nonlinear Dynamical Systems}.
World Scientific

\bibitem{fax2001graph}
Fax JA, Murray RM. 2001.
Graph {Laplacians} and stabilization of vehicle formations

\bibitem{Mirzaev2013LaplacianDynamics}
Mirzaev I, Gunawardena J. 2013.
Laplacian dynamics on general graphs.
\textit{Bulletin of Mathematical Biology} 75:2118--2149

\bibitem{dorfler2018electrical}
D{\"o}rfler F, Simpson-Porco JW, Bullo F. 2018.
Electrical networks and algebraic graph theory: Models, properties, and
  applications.
\textit{Proceedings of the IEEE} 106:977--1005

\bibitem{Schaft2017modeling}
van~der Schaft AJ. 2017.
Modeling of physical network systems.
\textit{Systems \& Control Letters} 101:21--27

\bibitem{Newman2006booknetwork}
Newman ME, Barab{\'a}si ALE, Watts DJ. 2006.
\textit{The Structure and Dynamics of Networks}.
Princeton University Press

\bibitem{Fagnani2017book}
Fagnani F, Frasca P. 2017.
\textit{Introduction to Averaging Dynamics over Networks}, vol. 472.
Springer

\bibitem{rantzer2015positive}
Rantzer A. 2015.
Scalable control of positive systems.
\textit{European Journal of Control} 24:72--80

\bibitem{rantzer2018tutorial}
Rantzer A, Valcher ME. 2018.
\textit{A tutorial on positive systems and large scale control}.
In \textit{2018 IEEE Conference on Decision and Control (CDC)}, pp.
  3686--3697. IEEE

\bibitem{ishizaki2015directed}
Ishizaki T, Kashima K, Girard A, Imura Ji, Chen L, Aihara K. 2015.
Clustered model reduction of positive directed networks.
\textit{Automatica} 59:238--247

\bibitem{Johnson1990Matrix}
Johnson CR. 1990.
\textit{Matrix Theory and Applications}, vol.~40.
American Mathematical Society

\bibitem{scardovi2008synchronization}
Scardovi L, Sepulchre R. 2008.
\textit{Synchronization in networks of identical linear systems}.
In \textit{Proceedings of the 47th IEEE Conference on Decision and Control},
  pp.  546--551. IEEE

\bibitem{consensus2010li}
Li Z, Duan Z, Chen G, Huang L. 2010.
Consensus of multiagent systems and synchronization of complex networks: a
  unified viewpoint.
\textit{IEEE Transactions on Circuits and Systems I: Regular Papers}
  57:213--224

\bibitem{JanWillems1976}
G.Willems J. 1976.
Realization of systems with internal passivity and symmetry constraints.
\textit{Journal of the Franklin Institute} 301:605--621

\bibitem{arcak2007passivity}
Arcak M. 2007.
Passivity as a design tool for group coordination.
\textit{IEEE Transactions on Automatic Control} 52:1380--1390

\bibitem{chopra2012output}
Chopra N. 2012.
Output synchronization on strongly connected graphs.
\textit{IEEE Transactions on Automatic Control} 57:2896--2901

\bibitem{ChengAuto2020Gramian}
Cheng X, Scherpen JMA. 2020.
Novel gramians for linear semistable systems.
\textit{Automatica} 115:108911

\bibitem{bhat1999lyapunov}
Bhat SP, Bernstein DS. 1999.
\textit{Lyapunov analysis of semistability}.
In \textit{Proceedings of the 1999 American Control Conference}, vol.~3, pp.
  1608--1612. IEEE

\bibitem{hui2009semistability}
Hui Q, Haddad WM, Bhat SP. 2009.
Semistability, finite-time stability, differential inclusions, and
  discontinuous dynamical systems having a continuum of equilibria.
\textit{IEEE Transactions on Automatic Control} 54:2465--2470

\bibitem{ChengCDC2017Digraph}
Cheng X, Scherpen JMA. 2017{\natexlab{a}}.
\textit{A New Controllability Gramian for Semistable Systems and its
  Application to Approximation of Directed Networks}.
In \textit{Proceedings of the 56th IEEE Conference on Decision and Control},
  pp.  3823--3828. Melbourne, Australia

\bibitem{ChengTAC2020Weight}
Cheng X, Yu L, Ren D, Scherpen JMA. 2020.
Reduced order modeling of diffusively coupled network systems: An optimal edge
  weighting approach.
\textit{arXiv preprint arXiv:2003.03559}

\bibitem{farina2011positive}
Farina L, Rinaldi S. 2011.
\textit{Positive Linear Systems: Theory and Applications}, vol.~50.
John Wiley \& Sons

\bibitem{ChengTAC2019Digraph}
Cheng X, Scherpen JMA. 2019.
Clustering-based model reduction of {Laplacian} dynamics with weakly connected
  topology.
\textit{IEEE Transactions on Automatic Control}

\bibitem{zhang2013AEP}
Zhang S, Cao M, Camlibel MK. 2013.
Upper and lower bounds for controllable subspaces of networks of diffusively
  coupled agents.
\textit{IEEE Transactions on Automatic control} 59:745--750

\bibitem{aguilar2017AEP}
Aguilar CO, Gharesifard B. 2017.
Almost equitable partitions and new necessary conditions for network
  controllability.
\textit{Automatica} 80:25--31

\bibitem{Petar2015CDC}
Mlinari\'c P, Grundel S, Benner P. 2015.
\textit{Efficient model order reduction for multi-agent systems using {QR}
  decomposition-based clustering}.
In \textit{Proceedings of the 54th IEEE Conference on Decision and Control},
  pp.  4794--4799

\bibitem{ChengECC2019weight}
Cheng X, Yu L, Scherpen JMA. 2019.
\textit{Reduced order modeling of linear consensus networks using weight
  assignments}.
In \textit{Proceedings of the 17th European Control Conference}, pp.
  2005--2010. Napoli, Italy

\bibitem{zelazo2013performance}
Zelazo D, Schuler S, Allg{\"o}wer F. 2013.
Performance and design of cycles in consensus networks.
\textit{Systems \& Control Letters} 62:85--96

\bibitem{leiter2017graph}
Leiter N, Zelazo D. 2017.
Graph-based model reduction of the controlled consensus protocol.
\textit{IFAC-PapersOnLine} 50:9456--9461

\bibitem{ChengECC2016}
Cheng X, Kawano Y, Scherpen JMA. 2016.
\textit{Graph structure-preserving model reduction of linear network systems}.
In \textit{Proceedings of the 15th European Control Conference}, pp.
  1970--1975. Aalborg, Denmark

\bibitem{ChengCDC2016Gramian}
Cheng X, Scherpen JMA. 2016.
\textit{Introducing network Gramians to undirected network systems for
  structure-preserving model reduction}.
In \textit{Proceedings of the 55th IEEE Conference on Decision and Control},
  pp.  5756--5761. Las Vegas, the USA

\bibitem{ChengACOM2018Power}
Cheng X, Scherpen JMA. 2018{\natexlab{a}}.
Clustering approach to model order reduction of power networks with distributed
  controllers.
\textit{Advances in Computational Mathematics} 44:1917--1939

\bibitem{kawano2019data}
Kawano Y, Besselink B, Scherpen JM, Cao M. 2020.
Data-driven model reduction of monotone systems by nonlinear dc gains.
\textit{IEEE Transactions on Automatic Control} 65:2094 -- 2106

\bibitem{zha2002spectral}
Zha H, He X, Ding C, Gu M, Simon HD. 2002.
\textit{Spectral relaxation for k-means clustering}.
In \textit{Advances in neural information processing systems}, pp.  1057--1064

\bibitem{UmarCDC2019Clustering}
Niazi MUB, Cheng X, Canudas~de Wit C, Scherpen JMA. 2019.
\textit{Structure-based clustering for model reduction of large-scale
  networks}.
In \textit{Proceedings of the 58th IEEE Conference on Decision and Control},
  pp.  5038--5043. Nice, France

\bibitem{aoki1968aggregation}
Aoki M. 1968.
Control of large-scale dynamic systems by aggregation.
\textit{IEEE Transactions on Automatic Control} 13:246--253

\bibitem{feliachi1987interconnected}
Feliachi A, Bhurtun C. 1987.
Model reduction of large-scale interconnected systems.
\textit{International Journal of Systems Science} 18:2249--2259

\bibitem{ishizaki2015clustered2o}
Ishizaki T, Imura JI. 2015.
Clustered model reduction of interconnected second-order systems.
\textit{Nonlinear Theory and Its Applications, IEICE} 6:26--37

\bibitem{ChengECC2020}
Cheng X, Necoara I, Lupu D. 2020.
\textit{A suboptimal $H_2$ clustering-based model reduction approach for linear
  network systems}.
In \textit{17th European Control Conference (ECC2020)}, pp.  1961--1966

\bibitem{lacoste2016FW}
Lacoste-Julien S. 2016.
Convergence rate of frank-wolfe for non-convex objectives.
\textit{arXiv preprint arXiv:1607.00345}

\bibitem{dinh2011combining}
Dinh QT, Gumussoy S, Michiels W, Diehl M. 2011.
Combining convex-concave decompositions and linearization approaches for
  solving {BMIs}, with application to static output feedback.
\textit{IEEE Transactions on Automatic Control} 57:1377--1390

\bibitem{kokotovic1976singular}
Kokotovic PV, O'Malley~Jr RE, Sannuti P. 1976.
Singular perturbations and order reduction in control theory—an overview.
\textit{Automatica} 12:123--132

\bibitem{anderson2011reduction}
Anderson J, Chang YC, Papachristodoulou A. 2011.
Model decomposition and reduction tools for large-scale networks in systems
  biology.
\textit{Automatica} 47:1165--1174

\bibitem{hancock2015simplified}
Hancock EJ, Stan GB, Arpino JA, Papachristodoulou A. 2015.
Simplified mechanistic models of gene regulation for analysis and design.
\textit{Journal of The Royal Society Interface} 12:20150312

\bibitem{chow1995inertial}
Chow JH, Galarza R, Accari P, Price WW. 1995.
Inertial and slow coherency aggregation algorithms for power system dynamic
  model reduction.
\textit{IEEE Transactions on Power Systems} 10:680--685

\bibitem{Biyik2008Areaaggregation}
B{\i}y{\i}k E, Arcak M. 2008.
Area aggregation and time-scale modeling for sparse nonlinear networks.
\textit{Systems \& Control Letters} 57:142--149

\bibitem{romeres2013novel}
Romeres D, D{\"o}rfler F, Bullo F. 2013.
\textit{Novel results on slow coherency in consensus and power networks}.
In \textit{Proceedings of the 2013 European Control Conference}, pp.  742--747.
  IEEE

\bibitem{Dorfler2013Kron}
D{\"o}rfler F, Bullo F. 2013.
Kron reduction of graphs with applications to electrical networks.
\textit{IEEE Transactions on Circuits and Systems I: Regular Papers}
  60:150--163

\bibitem{Monshizadeh2017Constant}
Monshizadeh N, De~Persis C, van~der Schaft AJ, Scherpen JMA. 2017.
A novel reduced model for electrical networks with constant power loads.
\textit{IEEE Transactions on Automatic Control}

\bibitem{deng2011optimal}
Deng K, Mehta PG, Meyn SP. 2011.
Optimal kullback-leibler aggregation via spectral theory of markov chains.
\textit{IEEE Transactions on Automatic Control} 56:2793--2808

\bibitem{Martin2018mergetocure}
Martin N, Frasca P, Canudas-De-Wit C. 2018.
Mergetocure: a new strategy to allocate cure in an epidemic over a grid-like
  network using a scale-free abstraction.
\textit{IFAC-PapersOnLine} 51:34--39

\bibitem{ChengIFAC2017}
Cheng X, Scherpen JMA. 2017{\natexlab{b}}.
Balanced truncation approach to linear network system model order reduction.
\textit{IFAC-PapersOnLine} 50:2451--2456

\bibitem{LanlinCDC2019second-order}
Yu L, Cheng X, Scherpen JMA. 2019{\natexlab{a}}.
\textit{$H_2$ sub-optimal model reduction for second-order network systems}.
In \textit{Proceedings of the 58th IEEE Conference on Decision and Control},
  pp.  5062--5067. Nice, France

\bibitem{LanlinCDC2019eigenvalue}
Yu L, Cheng X, Scherpen JMA. 2019{\natexlab{b}}.
\textit{Synchronization preserving model reduction of multi-agent network
  systems by eigenvalue assignments}.
In \textit{Proceedings of the 58th IEEE Conference on Decision and Control},
  pp.  7794--7799. Nice, France

\bibitem{Fujimoto2010BT}
Fujimoto K, Scherpen JMA. 2010.
Balanced realization and model order reduction for nonlinear systems based on
  singular value analysis.
\textit{SIAM Journal on Control and Optimization} 48:4591--4623

\bibitem{Scherpen11}
Scherpen JMA. 2011.
Balanced realization, balanced truncation and the {H}ankel operator, chap.~4.
The Control Handbook, Control Systems Advanced Methods, Eds. W. Levine, CRC
  Press, Taylor \& Francis Group,  1--24

\bibitem{ishizaki2013structured}
Ishizaki T, Sandberg H, Johansson KH, Kashima K, Imura Ji, Aihara K. 2013.
\textit{Structured model reduction of interconnected linear systems based on
  singular perturbation}.
In \textit{2013 American Control Conference}, pp.  5524--5529. IEEE

\bibitem{jaoude2017balanced}
Jaoude DA, Farhood M. 2017.
Balanced truncation model reduction of nonstationary systems interconnected
  over arbitrary graphs.
\textit{Automatica} 85:405--411

\bibitem{Ionescu2020H2}
Necoara I, Ionescu TC. 2020.
{$ H_2 $} model reduction of linear network systems by moment matching and
  optimization.
\textit{IEEE Transactions on Automatic Control}

\bibitem{ChengECC2018Lure}
Cheng X, Scherpen JMA. 2018{\natexlab{b}}.
\textit{Robust Synchronization Preserving Model Reduction of {Lur'e} Networks}.
In \textit{Proceedings of the 16th European Control Conference}, pp.
  2254--2259. Limassol, Cyprus

\bibitem{Scarciotti17}
Scarciotti G, Astolfi A. 2005.
\textit{Nonlinear Model Reduction by Moment Matching}.
NOW Publishers, Foundation and Trends in Systems and Control, Vol. 4, No. 3-4,
  224-409

\bibitem{forrow2018spectra}
Forrow A, Woodhouse FG, Dunkel J. 2018.
Functional control of network dynamics using designed {Laplacian} spectra.
\textit{Physical Review X} 8:041043

\bibitem{jongsma2017CycleRemoval}
Jongsma HJ, Trentelman HL, Camlibel KM. 2017.
Model reduction of networked multiagent systems by cycle removal.
\textit{IEEE Transactions on Automatic Control} 63:657--671

\end{thebibliography}

\end{document}